\documentclass[reqno,11pt]{amsart}

\usepackage{graphicx}         
\usepackage{amsfonts}
\usepackage{amssymb}
\usepackage[latin1]{inputenc}
\usepackage[all]{xy}

\newtheorem{thm}{Theorem}[section]
\newtheorem{cor}[thm]{Corollary}
\newtheorem{lemma}[thm]{Lemma}

\newtheorem{prop}[thm]{Proposition}
\theoremstyle{definition}
\newtheorem{defn}[thm]{Definition}
\theoremstyle{remark}
\newtheorem{rem}[thm]{Remark}
\newtheorem{remark}[thm]{Remark}
\newtheorem{example}[thm]{Example}

\renewcommand{\a}{\alpha}
\renewcommand{\b}{\beta}

\renewcommand{\l}{\lambda}
\newcommand{\s}{\sigma}

\newcommand{\can}[1]{\left\langle#1\right\rangle}

\newcommand{\PB}{\{\cdot\,,\cdot\}}

\newcommand{\pb}[1]{\left\{#1\right\}}

\newcommand{\lb}[1]{\[#1\]}
\newcommand{\LB}{[\cdot\,,\cdot]}

\newcommand{\inn}[2]{\left\langle#1\,\vert\,#2\right\rangle}

\renewcommand{\[}{\left[}
\renewcommand{\]}{\right]}

\newcommand{\set}[1]{\left\{#1\right\}}

\newcommand{\cE}{\mathcal E} 
\newcommand{\cL}{\mathcal L}
\renewcommand{\L}{\mathcal L}
\newcommand{\cF}{\mathcal F}

\newcommand{\cV}{\mathcal V} 
\newcommand{\cW}{\mathcal W} 

\newcommand{\fA}{\mathfrak A} 
\newcommand{\fg}{\mathfrak g} 
\newcommand{\fu}{\mathfrak u} 
\newcommand{\fF}{\mathfrak F} 
\newcommand{\fG}{\mathfrak G} 
\newcommand{\fH}{\mathfrak H} 
\newcommand{\fT}{\mathfrak T} 

\newcommand{\bbP}{\mathbb P}

\newcommand{\bbR}{\mathbb R}
\newcommand{\bbT}{\mathbb T}

\newcommand{\bbZ}{\mathbb Z}

\newcommand{\F}{\mathbf F}
\newcommand{\G}{\mathbf G}

\newcommand{\X}{ X}

\newcommand{\Vect}{{\mathfrak{X}}}
\newcommand{\Ker}{\mathop{\rm Ker}\nolimits}
\newcommand{\Tr}{\mathop{\rm Trace}\nolimits}
\newcommand{\Ad}{\mathop{\rm Ad}\nolimits}

\newcommand{\HOM}{\mathop{\mathcal{H}om}\nolimits}
\renewcommand{\Im}{\mathop{\rm Im}\nolimits}
\newcommand{\leqs}{\leqslant}
\newcommand{\geqs}{\geqslant}
\newcommand{\pp}[2]{\frac{\partial#1}{\partial#2}}
\newcommand{\p}{\partial}
\newcommand{\we}{\wedge}
\newcommand{\Rk}{\hbox{rank\,}}
\newcommand{\diff}{{\rm d }}
\newcommand{\CASM}{{\mathcal Cas}^{M}_B}
\newcommand{\Cas}{{\mathcal Cas}}

\DeclareMathOperator{\Obs}{Obs}

\newif\ifprivate
\privatefalse

 \numberwithin{equation}{section}

\def\???{\ifprivate {\bf {???}} \marginpar{{\Huge {\bf ?}}}\else \fi}
\numberwithin{equation}{section}

\begin{document}  

\nocite{*} 

%
\title[Global Action-Angle Variables]{Global Action-Angle Variables for Non-Commutative Integrable Systems}

\author[Fernandes]{Rui L.~Fernandes} \address{Rui L.~Fernandes, Department of Mathematics, University of Illinois at 
        Urbana-Champaign, 1409 W. Green Street, Urbana, IL 61801, USA }\email{ruiloja@illinois.edu}

\author[Laurent-Gengoux]{Camille Laurent-Gengoux}\address{Camille Laurent-Gengoux, Laboratoire de Math\'ematiques, UMR 7122 du CNRS,
  Université de Metz, Ile du Saulcy, F-57045 Metz Cedex 1, France}\email{camille.laurentgengoux@univ-metz.fr}

\author[Vanhaecke]{Pol Vanhaecke} \address{Pol Vanhaecke, Laboratoire de Math\'ematiques et Applications, 
         UMR 7348 du CNRS, Universit\'e de Poitiers, Boulevard Marie et Pierre Curie, BP 30179, 86962 Futuroscope
         Chasseneuil Cedex, France}\email{pol.vanhaecke@math.univ-poitiers.fr}

\date{\today}



\begin{abstract}

In this paper we analyze the obstructions to the existence of global action-angle variables for regular
non-commutative integrable systems (NCI systems) on Poisson manifolds.  In contrast with local action-angle
variables, which exist as soon as the fibers of the momentum map of such an integrable system are compact, global
action-angle variables rarely exist. This fact was first observed and analyzed by Duistermaat in the case of
Liouville integrable systems on symplectic manifolds and later by Dazord-Delzant in the case of non-commutative
integrable systems on symplectic manifolds. In our more general case where phase space is an arbitrary Poisson
manifold, there are more obstructions, as we will show both abstractly and on concrete examples.  Our approach
makes use of a few new features which we introduce: the action bundle and the action lattice bundle of the NCI
system (these bundles are canonically defined) and three foliations (the action, angle and transverse foliation),
whose existence is also subject to obstructions, often of a cohomological nature.

\end{abstract}

\maketitle

\tableofcontents

\section{Introduction}
The notion of a Liouville integrable system on a symplectic manifold \cite[Ch.~10]{arnold} has two natural
generalizations, namely the notion of a Liouville integrable system on a Poisson manifold
\cite[Ch.\ 4]{adlermoerbekevanhaecke2004} and the notion of a non-commutative integrable system on a symplectic
manifold \cite{bolsinov,FA,AF,Nekhroshev}. These two concepts were merged in \cite{LMV}, where the notion of a
non-commutative integrable system on a Poisson manifold was introduced.

A \textbf{non-commutative integrable system} (NCI system) on an $n$-dimen\-sional Poisson manifold $(M,\Pi)$ is a
family $f_1,\dots,f_s$ of smooth functions on~$M$, such that the first $n-s$ functions are in involution (Poisson
commute) with every function in the family:
$$ 
  \{f_i,f_j \} =0 , \hbox{ for } 1\leqslant i\leqslant n-s,\ 1\leqslant j\leqslant s\;,
$$
and satisfy an independence condition which will be stated below. The number $r:=n-s$ is called the \textbf{rank}
of the NCI system.  The classical case of Liouville integrable systems on a symplectic manifold corresponds to the
case where $r=s=n/2$, while the case of superintegrable systems (on a symplectic or Poisson manifold) corresponds
to $r=1$; for other NCI systems, $r$ can be any integer satisfying $2\leqslant 2r\leqslant n$.

One usually thinks of an NCI system on an $n$-dimensional phase space~$M$ as a Hamiltonian dynamical system $\X_h$
on $M$, associated with some function~$h$, admitting the functions $h=f_1,f_2,\dots,f_s$, as \textbf{first
integrals}, i.e., $X_hf_i=\pb{f_i,h}=0$ for $1\leqs i\leqs s$. Then the above definition of an NCI system can be
understood as follows: (i) one can first reduce the dynamics of $\X_h$ to a generic common level set of all the
first integrals $f_1,\dots,f_s$, thereby reducing the dimension of the phase space by $s$; (ii) since the first
integrals $f_1,\dots,f_r$ are in involution with all the above first integrals, the flows of the Hamiltonian vector
fields $\X_{f_1},\dots,\X_{f_r}$ define a local $\bbR^r$-action which preserves this common level set, so one can
further reduce the dimension of the system by $r$ by passing to the quotient space of the level set by the
action. Altogether, one can reduce the dimension by $r+s=n$, the dimension of the phase space, which justifies the
name ``integrable''. To be precise, the above dimension count is correct only if we assume independence of the
first integrals:
\begin{equation}\label{eq:intro_indep1}
  \diff f_1\wedge \cdots \wedge \diff f_s\not=0\;,
\end{equation}
as well as of the Hamiltonian vector fields generating the local $\bbR^r$-action:
\begin{equation}\label{eq:intro_indep2}
  \X_{f_1}\wedge \cdots \wedge \X_{f_r}\not=0\;.
\end{equation}
The latter condition does, in general, not follow from the former condition because the Poisson tensor may have a
non-trivial kernel. We will deal in this paper solely with \textbf{regular NCI systems}, i.e., NCI systems such
that conditions (\ref{eq:intro_indep1}) and (\ref{eq:intro_indep2}) hold at every point of $M$. The study of
singularities of NCI systems (points where at least one of the above conditions fails) is a very important and
interesting topic, which we defer to future works.


Examples of NCI systems include, besides Liouville integrable systems, many classical systems such as the motion in
a central force field, the Kepler problem, the Euler-Poinsot top and the Gelfand-Cetlin system. Each one of these
systems has singularities, but by removing some appropriate closed subset which contains them, we obtain a regular
NCI system to which the theory developed here applies.

We assemble the first integrals of an NCI system in a single map $\F=(f_1,\dots,f_s):M\to \bbR^s$, which we call
the \textbf{momentum map} of the NCI system. Notice that $\F$ is submersive when the NCI system is regular. 
The first important, non-trivial, fact about NCI systems on Poisson manifolds is the action-angle theorem, which
was proved in full generality in \cite{LMV}. We state it here for regular NCI systems for which the fibers of its
momentum map are compact and connected.

\begin{thm}[Existence of local action-angle variables]
  Let $(M,\Pi,\F)$ be a regular NCI system of dimension $n$ and rank $r=n-s$ with compact connected fibers. For any
  $b$ in the image of $\F$, there exists an open neighborhood~$U$ of $b$ in $\bbR^{s}$, an open neighborhood $V$ of
  $\F^{-1}(b)$ in $M$ and an open embedding $\Psi:V\to T^*\bbT^r\times\bbR^{s-r}$ such that the following diagram
  is commutative:
  $$
    \xymatrix{\F^{-1}(U)\supset V\ar[r]^{\Psi}\ar[d]_{\F_{\vert V}}& T^*\bbT^r\times\bbR^{s-r}\ar[d]\\
            U\ar[r]&\bbR^{s}}
  $$
  Moreover, $\Psi$ is a Poisson map if we consider on $T^*\bbT^r\times\bbR^{s-r}$ the product of the canonical
  symplectic structure on $T^*\bbT^r$ with an appropriate Poisson structure on an open subset of $\bbR^{s-r}$.
\end{thm}

The above theorem is semi-local in the sense that it describes such NCI systems in the neighborhood of a connected
component of a fiber (of the map~$\F$; such a component is an $r$-dimensional torus $\bbT^r$, just like in the
classical Liouville theorem), rather than in the neighborhood of a point. In terms of the natural coordinates
$(\theta_i,p_i,z_j)$ on $ \bbT^r \times \bbR^r\times\bbR^{n-2r}\simeq T^*\bbT^r\times\bbR^{s-r}$ the Poisson
structure on $M$ takes the following form:
$$ \Pi= \sum_{i=1}^r\frac{\partial }{\partial \theta_i}\wedge \frac{\partial }{\partial p_i}+\sum_{1\leqslant
  j<k\leqslant s-r}c_{jk}(z) \frac{\partial }{\partial z_j}\wedge \frac{\partial }{\partial z_k}\;,
$$ 
where the second sum is absent in the case of Liouville integrable systems on regular Poisson manifolds, such as
symplectic manifolds. The variables $\theta_i,p_i$ and $z_j$, in that order, are called \textbf{angle, action} and
\textbf{transverse variables} (or \textbf{coordinates}; it is understood that the $\theta_i$ are $S^1$-valued).

According to the theorem, the phase space of a regular NCI system (with compact connected fibers) can be covered
with charts equipped with action-angle-transverse variables. Of course, these local variables are highly
non-unique. Therefore, the question asking whether for a given NCI system these local variables can be glued to
yield global variables is a non-trivial one. The main focus in this paper is to describe the different obstructions
for this passage from local to global. As an intermediate step, we will also consider the obstructions to the
existence of action, angle and transverse foliations, which are weaker than the obstructions for the existence of
the corresponding variables, but are in general easier to compute. For each of these obstructions, we will prove
their non-triviality in some concrete examples.

We now give an outline of the paper and describe the main results.

In Section \ref{sec:NCI} we recall the notion of a non-commutative integrable system on a Poisson manifold, which
we reformulate in geometrical terms (in terms of a foliation) and we initiate the study of the Poisson geometry of
such a system. The upshot is that we view regular NCI systems as Poisson maps $(M,\Pi)\stackrel{\phi}{\to
}(B,\pi)$, whose fibers define a rank $r$ foliation with compact leaves.

A key novelty which is introduced in Section~\ref{sec:action} is the \textbf{action bundle} $E$, which is a vector
bundle of rank $r$ on $B$ and whose sections generate, upon using the Poisson structure $\Pi$, the action vector
fields, i.e., the commuting, integrable vector fields which are tangent to the fibers of the momentum map
(Proposition \ref{prp:ranks} and Lemma \ref{lma:fund_vfs}). When the fibers of the momentum map $\phi$ are compact
and connected, $E$ contains a lattice bundle $L_B\to B$, the \textbf{action lattice bundle}, whose sections
generate periodic vector fields of period~1; it implies that $M$ is a torus bundle over $B$ (see Section
\ref{par:lattice}). A set of action variables of the NCI system is a collection of $r$ functions on $B$ which
define a global trivialization of the action lattice bundle (making the torus bundle $\phi:M\to B$ into a principal
$\bbT^r$-bundle); the obstruction to their existence lies in $H^1(B,\CASM)$, where $\CASM$ is the sheaf of
functions on $B$ who pull back to Casimir functions on~$M$ (Theorem \ref{thm:action_var_obs}). Action variables
define a (transversely integral affine) foliation on $B$, which leads to the notion of an \textbf{action
foliation}.  When the action lattice bundle admits a trivialization on $B$, it defines a cohomology class in
$H^1(B,\CASM/\bbR)$, whose nullity is equivalent to the existence of an action foliation. This class is, of course,
closely related to $H^1(B,\CASM)$, which is decisive for the existence of action variables (Proposition
\ref{cor:actionfoliation}).

The existence of angle variables is discussed in Section \ref{sec:angle}. Interestingly, they can be defined in
terms of the action lattice bundle, hence their (global) existence can be studied independently of the existence of
action variables, or of a choice of such variables. We show that global angle variables exist if and only if the
action lattice bundle is trivial and the momentum map $\phi:M\to B$ admits a coisotropic section (Theorem
\ref{thm:angle_existence}). The latter condition is in an essential way non-linear, hence does not lead to a
cohomological obstruction class, as in the case of action variables. However, a set of angle variables defines a
pair of foliations, an \textbf{angle foliation} (which is a foliation of $M$) and a \textbf{transverse foliation}
(which is a foliation of $B$, transverse to every action foliation, see Propositions
\ref{prop:angles_distribution_quotient} and \ref{prop:angle_foliation_distribution_quotient}). The obstructions to
the existence of such a pair of foliations then leads to obstructions of the existence of global angle variables,
which are weaker than the existence of a coisotropic section, but easier to compute explicitly. We finish Section
\ref{sec:angle} with a theorem which gives an explicit description of every NCI system for which action-angle
variables do exist, under the assumption that all leaves of the action and the transverse foliation intersect in a
unique point (Theorem \ref{thm:aa_when_exists}): in terms of angle variables $\theta_i$ and action variables $p_i$,
the Poisson structure on its phase space $M$ then takes the canonical form
\begin{equation*}
    \Pi=\sum_{i=1}^r \frac{\partial}{\partial \theta_i} \wedge \frac{\partial}{\partial p_i} + \pi_{\vert_A}\;,
\end{equation*}
where $A$ is any leaf of the action foliation (which turns out to be a Poisson submanifold of $B$).

Section \ref{section:examples} is devoted to the study of several examples. They include articifially constructed
mathematical examples which illustrate the non-triviality of the obstructions that are discussed in the paper, as
well as examples coming from classical mechanics, which turn out to exhibit a large spectrum of phenomena which
have a definite impact on the global geometry of NCI systems.

%
\bigskip

\begin{center}
{\bf Conventions}
\end{center}

In this paper, all manifolds and objects considered on them are real and smooth. When $\Pi$ is a Poisson structure
on a manifold $M$, we write $\pb{f,g}$ for $\Pi(\diff f,\diff g)$ and we denote the Hamiltonian vector field
associated to $h\in C^\infty(M)$ by~$\X_h$. The vector bundle map induced by $\Pi$ is denoted by
$\Pi^\sharp:T^*M\to TM$. Our sign convention is that $\X_h(g)=\diff g(\X_h)=\pb{g,h}$ for $g\in C^\infty(M)$ and
$\Pi^\sharp(\diff h)=-\X_h$.  For a foliation $\fF$ on a manifold $M$ the tangent space to~$\fF$ at~$m$ is denoted
by $T_m\fF$, while its annihilator is denoted by $(T_m\fF)^\circ$. It leads to subbundles $T\fF$ of $TM$ and
$(T\fF)^\circ$ of $T^*M$. For a vector bundle~$E$ over $M$, the module of (smooth) sections of $E$ is denoted by
$\Gamma(E)$. We denote by $\Omega^k(M)$ (respectively by $\Vect^k(M)$) the module $\Gamma(\we^kT^*M)$ of $k$-forms
(respectively the module $\Gamma(\we^kTM)$ of $k$-vector fields) on $M$. For $\omega\in\Omega^k(M)$ we denote by
$\omega_m$ or $\omega\vert_m$ its value at $m\in M$ and similarly for elements of $\Vect^k(M)$. For a vector
field~$\cV$ on $M$, we denote by $\L_\cV$ the Lie derivative with respect to $\cV$ of elements of $\Omega^k(M)$ or
of $\Vect^k(M)$. The $r$-dimensional torus $(\bbR/\bbZ)^r$ is denoted by $\bbT^r$.

%
\section{Non-commutative integrable systems on Poisson manifolds}\label{sec:NCI}

\subsection{NCI systems} We first  recall from \cite{LMV}  the main notion relevant to this paper.
\begin{defn}\label{def:non-com}
Let $(M,\Pi)$ be a Poisson manifold of dimension $n$. Let $\F=(f_1,\dots,f_s)$ be an $s$-tuple of functions on $M$,
where $2s\geqs n$ and set $r:=n-s$. Suppose the following:
\begin{enumerate}
  \item[(1)] The functions $f_1,\dots,f_r$ are in involution with the functions $f_1,\dots,f_s$:
  $$ \{f_i,f_j\}=0,\qquad (1\leqs i\leqs r \hbox{ and } 1\leqs j \leqs s)\;;$$
  \item[(2)] For $m$ in a dense open subset of $M$:
  $$ \diff_m f_1\wedge \dots \wedge \diff_m f_s\ne 0 \quad \hbox{and} \quad \X_{f_1}|_m\wedge \dots \wedge
    \X_{f_r}|_m\ne 0\;.  $$
\end{enumerate}
Then the triplet $(M,\Pi,\F)$ is called a \textbf{non-commutative integrable system} ({NCI system}) of
\textbf{rank} $r$ and $\F$, viewed as a map $\F:M\to\bbR^s$, is called its \textbf{momentum map}.
\end{defn}

The classical case of a \textbf{Liouville integrable system} corresponds to the particular case where $r$ is half
the (maximal) rank of $\Pi$; this implies that \emph{all} the functions $f_1,\dots,f_s$ are pairwise in involution,
$$ 
  \{f_i,f_j\}=0\qquad (1\leqs i,j \leqs s)\;.
$$

A point $m\in M$ where the two conditions in (2) hold is called a \textbf{regular point} of the NCI system, the
other points are called \textbf{singular points} of the NCI system. When all points of $M$ are regular one speaks
of a \textbf{regular NCI system}. We will mainly study regular NCI systems, though we will see in Section
\ref{section:examples} that singular points are present in basically all the examples; we will then be led to
restricting the Poisson manifold underlying the NCI system to an appropriate open subset, on which the NCI system
restricts to a regular NCI system.

We start with an example from classical mechanics (see \cite[Ch.\ 4.48]{whittaker}).
\begin{example}
Consider a particle of mass $m$ in $\bbR^3$ which is subject to a central force, derived from a potential function
$V=V(r)$ which depends only on the distance $r$ from the origin of $\bbR^3$. The Hamiltonian which describes the
total energy of the particle is given by
\begin{equation*}
  H=\frac1{2m}\sum_{i=1}^3p_i^2+V(r)\;,
\end{equation*}%
where $r^2=\sum_{i=1}^3q_i^2$ and where $(q_1,q_2,q_3)$ and $(p_1,p_2,p_3)$ respectively stand for the position
coordinates and for the corresponding momenta of the particle. The Poisson structure is the canonical structure on
$T^*\bbR^3\simeq\bbR^6$, to wit
$$
  \Pi=\sum_{i=1}^3\pp{}{q_i}\we\pp{}{p_i}\;.
$$
The Hamiltonian vector field $\X_H$ whose integral curves describe the motion of the particle is given by
\begin{equation}\label{eq:motion_central}
  \dot q_i=\pp H{p_i}=\frac{p_i}m\;,\qquad
  \dot p_i=-\pp H{q_i}=-q_i\frac{V'(r)}r\;.
\end{equation}%
Consider the three linear momenta $\mu_{ij}:=q_ip_j-q_jp_i$, where $1\leqslant i<j\leqslant 3$. It follows at once
from (\ref{eq:motion_central}) that $\dot\mu_{ij}=0$, so that each of these momenta is a constant of motion, and so
$L:=\mu_{12}^2+\mu_{13}^2+\mu_{23}^2$ is also a constant of motion; moreover, the latter has the virtue of being in
involution with all the linear momenta $\mu_{ij}$. Letting $\F:=(H,L,\mu_{12},\mu_{23})$ it follows that
$(T^*\bbR^3,\Pi,\F)$ is an NCI system of rank 2 with momentum map $\F$.
\end{example}
Next, we give a family of examples of regular NCI systems which are important for the theory which will be
developed in this paper, because they to provide local models for any regular NCI system (see
Proposition~\ref{prp:NCI_local} below).
\begin{example}\label{ex:can_model}
Let $M:=\bbR^{2r}\times\bbR^{s-r}$ with coordinates $(q_i,p_i,z_j)$ be equipped with a Poisson structure $\Pi$ of
the form:
$$ 
  \Pi= \sum_{i=1}^r\frac{\partial }{\partial q_i}\wedge \frac{\partial }{\partial p_i}+\pi\;,
$$ 
where $\pi$ is any Poisson structure on $\bbR^{s-r}$,
\begin{equation}\label{eq:small_poisson}
  \pi=\sum_{1\leqslant j<k\leqslant s-r}c_{jk}(z) \frac{\partial }{\partial z_j}\wedge \frac{\partial }{\partial
    z_k}\;.
\end{equation}
Letting $\F:=(p_1,\dots,p_r,z_1,\dots,z_{s-r})$ it is clear that $(M,\Pi,\F)$ is a regular NCI system of rank $r$
with momentum map $\F$. It is a Liouville integrable system if and only $\pi=0$ (equivalently, all functions
$c_{ij}$ are zero).
\end{example}

A slight modification of this example yields a family of examples of regular NCI systems with compact fibers, which
are semi-local models for regular NCI systems with compact fibers (see Theorem \ref{thm:aa_semi_local} below).

\begin{example}\label{ex:can_model_semi}
Let $M=T^*\bbT^r\times\bbR^{s-r}\simeq \bbT^r \times \bbR^r\times \bbR^{s-r}$ with coordinates $(\theta_i,p_i,z_j)$
be equipped with a Poisson structure $\Pi$ of the form:
$$ 
  \Pi= \sum_{i=1}^r\frac{\partial }{\partial \theta_i}\wedge \frac{\partial }{\partial p_i}+\pi\;,
$$ 
where $\pi$ is any Poisson structure on $\bbR^{s-r}$, as in (\ref{eq:small_poisson}). Letting $\F:=(p_1,\dots,$
$p_r,z_1,\dots,z_{s-r})$ we have as above that $(M,\Pi,\F)$ is a regular NCI system of rank $r$ with momentum map
$\F$.
\end{example}

\subsection{Abstract NCI systems}
To a regular NCI system $(M,\Pi,\F)$ one naturally associates an $r$-dimen\-sional foliation of $M$: by the
regularity assumption, $\F:M\to \bbR^s$ is a submersion onto some open subset $B\subset \bbR^s$, so that the
connected components of the fibers of $\F$, which are $r$-dimensional, are the leaves of a foliation $\fF$ of
$M$. In the case of Example \ref{ex:can_model} (resp.\ Example \ref{ex:can_model_semi}), these leaves are
$r$-dimensional affine spaces $\bbR^r$ (resp.\ $r$-dimensional tori~$\bbT^r$).

In the following proposition we rewrite the key elements of the definition of a regular NCI system in terms of the
foliation which is associated to it. Before doing this, let us recall that a (locally defined) function which is
constant on the leaves of a foliation $\fF$ is called a (local) \textbf{first integral} of $\fF$. These functions
are characterized by the property that they are annihilated by any set of vector fields which generate the tangent
bundle $T\fF$ to $\fF$. In Example \ref{ex:can_model} (resp.\ Example \ref{ex:can_model_semi}), the first
integrals of the foliation defined by~$\F$ are the functions on $M$ which are independent of $q_1,\dots,q_r$
(resp.\ of $\theta_1,\dots,\theta_r$).

\begin{prop}
\label{prop:NCIS}
Let $(M,\Pi,\F)$ be a regular NCI system of dimension $n$ and rank~$r$ and let $\fF$ denote the foliation whose
leaves are the connected components of the fibers of its momentum map $\F:M\to B$. Then $T\fF$ is spanned by
Hamiltonian vector fields associated to first integrals of $\fF$, i.e., for each $m\in M$ there exist local first
integrals of $\fF$, namely $f_1,\dots,f_r$, whose Hamiltonian vector fields span $T_{m'}\fF$, for $m'$ in a
neighborhood of $m$ in $M$. In particular, every leaf of $\fF$ is contained in a symplectic leaf of $\Pi$.
\end{prop}

\begin{proof}
Item (1) in Definition \ref{def:non-com} implies that the Hamiltonian vector fields $\X_{f_1},\dots,\X_{f_r}$ are
tangent to the fibers of $\F:M\to B$ (i.e., to the leaves of $\fF$), while item (2) implies that they actually span
the tangent spaces to these fibers at every regular point, i.e., at every point (since it is assumed that the NCI
system is regular). This shows that $T\fF$ is spanned by the Hamiltonian vector fields associated to the first
integrals $f_1,\dots f_r$ of $\fF$. As a consequence, every leaf of $\fF$ is contained in a symplectic leaf of
$\Pi$.
\end{proof}

The above proposition leads to the following more abstract notion of an NCI system and of morphisms between such
systems:

\begin{defn}\label{def:ANCIS}
Let $ (M,\Pi)$ be a Poisson manifold. An \textbf{abstract non-commutative integrable system} (abstract NCI system)
of rank $r$ is an $r$-dimensional foliation $\fF$ of $M$, whose tangent bundle $T\fF$ is spanned by Hamiltonian
vector fields associated to (local) first integrals of $\fF$.

A \textbf{morphism} between two abstract NCI systems $(M,\Pi,\fF)$ and $(N,\Theta,\fG)$ is a Poisson map $\phi:M\to
N$ which is transverse to $\fG$ and such that $\phi^{*}(\fG)=\fF$.
\end{defn}


\begin{example} 
A \textbf{Lagrangian foliation} of a Poisson manifold $(M,\Pi)$ is a foliation $\fF$ of $M$ for which
$T\fF=\Pi^\sharp(T\fF)^\circ$. In the terminology of Definition~\ref{def:isotropic} and Example
\ref{exa:coisotropic} this amounts to saying that $\fF$ is both isotropic and coisotropic; it implies that $\Pi$ is
regular, of rank twice the dimension of $\fF$. For a Lagrangian foliation $\fF$, the Hamiltonian vector fields
associated to all its first integrals are both tangent to $T\fF$ and span~$T\fF$. In particular, $(M,\Pi,\fF)$ is
an abstract NCI system. It is the abstract version of a Liouville integrable system (on a regular Poisson
manifold).
\end{example}

\begin{example}
Let $(M,\Pi)$ be a Poisson manifold. Any nowhere vanishing Hamiltonian vector field $\X_h$ defines a foliation
$\fF$, making $(M,\Pi,\fF)$ into an abstract NCI system of rank 1. In this case, the first integrals of $\fF$ are
precisely the first integrals of $\X_h$.
\end{example}

In view of Proposition \ref{prop:NCIS}, if $(M,\Pi,\F)$ is a regular NCI system and $\fF$ its associated foliation,
then $(M,\Pi,\fF)$ is an abstract NCI system. We show that, \emph{locally}, the converse is also true. We do this
by showing that locally every regular NCI system is isomorphic to one of the systems described in Example
\ref{ex:can_model}.
\begin{prop}\label{prp:NCI_local}
  Let $(M,\Pi,\fF)$ be an abstract NCI system of dimension~$n$ and rank $r$. Let $m$ be an arbitrary point of
  $M$. There exist on a neighborhood~$U$ of $m$ coordinates $q_1,\dots,q_r,p_1,\dots,p_r,z_1,\dots,z_{n-2r}$ such
  that the foliation~$\fF$ is defined on $U$ by the functions $p_1,\dots,p_r,z_1,\dots,z_{n-2r}$ and such that
  $\Pi$ is given, on $U$, by
  \begin{equation}\label{eq:Pi_local_nci}
    \Pi= \sum_{i=1}^r\frac{\partial }{\partial q_i}\wedge \frac{\partial }{\partial p_i}+ \sum_{1\leqslant
      j<k\leqslant n-2r}c_{jk}(z) \frac{\partial }{\partial z_j}\wedge \frac{\partial }{\partial z_k}\;,
  \end{equation}
  where the functions $c_{jk}$ are independent of $q_1,\dots,q_r,p_1,\dots,p_r$. In particular, setting
  $\F:=(p_1,\dots,p_r,z_1,\dots,z_{n-2r})$ we have that $(U,\Pi_{\vert U},\F)$ is a regular NCI system of rank $r$.
\end{prop}
\begin{proof}
The proof is a direct application of the Carathéodory-Jacobi-Lie theorem for Poisson manifolds (see
\cite[Sect.~2]{LMV} for a proof). This theorem says that if $(M,\Pi)$ is any Poisson manifold of dimension $n$ on
which $r$ functions $p_1,\dots,p_r$ are given, which are pairwise in involution and have independent Hamiltonian
vector fields at some point $m\in M$, then these functions can be extended to a coordinate system
$q_1,\dots,q_r,p_1,\dots,p_r,z_1,\dots,z_{n-2r}$ on a neighborhood $U$ of $m$, such that $\Pi$ takes on $U$ the
form (\ref{eq:Pi_local_nci}). In order to apply this theorem in the present case, we take any point $m$ of $M$ and
we choose as functions $p_1,\dots,p_r$ local first integrals of $\fF$ whose Hamiltonian vector fields generate
$T\fF$ in a neighborhood of $m$. These $r$ functions are in involution so the theorem applies. Notice that in view
of (\ref{eq:Pi_local_nci}) the tangent space to $\fF$ is spanned by the vector fields $\p/\p q_1,\dots,\p/\p q_r$,
so the first integrals of $\fF$ are the functions which are independent of $q_1,\dots,q_r$ and~$\fF$ is locally
defined by the functions $p_1,\dots,p_r,z_1,\dots,z_{n-2r}$.
\end{proof}
In order to give another example of an abstract NCI system, we need a result which is interesting in its own right.
\begin{cor}
\label{cor:ham_tgt}
Let $(M,\Pi,\fF)$ be an abstract NCI system of dimension~$n$ and rank~$r$.  If $\cV$ is a Hamiltonian vector field
which is tangent to the fibers of $\fF$, then every Hamiltonian of $\cV$ is a first integral of $\fF$;
\end{cor}
\begin{proof}
The proof follows at once from Proposition \ref{prp:NCI_local}. We give a direct proof. Let $m$ be an arbitrary
point of $M$. In view of the definition of an abstract NCI systems, there exist on a neighborhood $U$ of $m$ first
integrals $f_1,\dots,f_r$ of $\fF$ whose Hamiltonian vector fields span $T\fF$ (on $U$). Thus, a function $f$ on
$U$ is a first integral of~$\fF$ if and only if $\X_{f_i}(f)=0$, for $i=1,\dots,r$. Suppose that $h$ is a function
on $M$ whose Hamiltonian vector field $\cV:=\X_h$ is tangent to $\fF$. Then
\begin{equation*}
  \X_{f_i}(h)=\pb{h,f_i}=-\X_h(f_i)=-\cV(f_i)=0\;,  
\end{equation*}%
so that $h$ is a first integral of $\fF$.  
\end{proof}
\begin{example}
Let $G\times M\to M$ be a Hamiltonian action of a Lie group~$G$ (with Lie algebra $\fg$) on a Poisson manifold
$(M,\Pi)$. Recall that this means that there exists a Lie algebra homomorphism $\mu:(\fg,\LB)\to (C^\infty(M),\PB)$
such that for every $x\in\fg$, the function $\mu(x)$ is a Hamiltonian for the fundamental vector field $\underline
x$ associated to $x$. We assume that the isotropy groups of the action have constant dimension, so that the orbits
are the leaves of a foliation $\fF$. We claim that the following conditions are equivalent:
\begin{enumerate}
  \item[(i)] $(M,\Pi,\fF)$ is an  abstract NCI system;
  \item[(ii)] For every $x\in\fg$, the function $\mu(x)$ is a first integral of $\fF$;
  \item[(iii)] $\mu([\fg,\fg])=0$.
\end{enumerate}
The implication (i) $\Rightarrow$ (ii) follows from Corollary \ref{cor:ham_tgt}, applied to the Hamiltonian
$\mu(x)$ of $\underline x$. Conversely, when (ii) holds $T\fF$ is spanned by the Hamiltonian vector fields
associated to certain first integrals of $\fF$, namely the functions $\mu(x)$ with $x\in\fg$, so $(M,\Pi,\fF$) is an
abstract NCI system. For $x,y\in\fg$ we have that
\begin{equation*}
  \underline y(\mu(x))=\pb{\mu(x),\mu(y)}=\mu([x,y])\;,
\end{equation*}%
from which the equivalence of (ii) and (iii) follows at once. Notice that (iii) is trivially satisfied when $\fg$
is abelian. Moreover, when the action is locally free, (iii) is equivalent to $[\fg,\fg]=0$, i.e., to $\fg$ being
abelian.
\end{example}

\subsection{Poisson complete isotropic foliations} 

The foliation of an abstract NCI system has two main features, which we first define and illustrate with some basic
examples.

\begin{defn}\label{def:isotropic}
  Let $(M,\Pi)$ be a Poisson manifold and suppose that $\fF$ is a foliation of $M$. 
    \begin{enumerate}
      \item[(1)] We say that $\fF$ is \textbf{Poisson complete} if the Poisson bracket of two (local) first
        integrals of $\fF$ is a (local) first integral of~$\fF$;
\item[(2)] We say that $\fF$ is \textbf{isotropic} if $T\fF\subset \Pi^\sharp(T\fF)^\circ$. 
\end{enumerate}
\end{defn}
\begin{example}
Let $(M,\omega)$ be a symplectic manifold and let $\Pi:=\omega^{-1}$ denote the Poisson structure corresponding to
$\omega$. Suppose that there exists a nowhere vanishing 1-form $\alpha$ on $M$. Then the corresponding vector
field~$\Pi^\sharp(\alpha)$ defines a foliation $\fF$ which is isotropic, since $\Pi^\sharp(\alpha)$ generates
$T\fF$ in every point of $M$; also, $\a\in(T\fF)^\circ$. If $\alpha\wedge\diff\alpha\neq0$ then this foliation is
not Poisson complete. Indeed, $\Pi^\sharp(T\fF)^\circ$, which is the symplectic orthogonal distribution to $\fF$,
coincides with $\Ker\alpha$, which is integrable if and only if $\alpha\wedge\diff\alpha=0$; but, as we will see in
Proposition~\ref{prop:Poisson foliations} below, if $\fF$ is Poisson complete then the distribution
$\Pi^\sharp(T\fF)^\circ$ is integrable. In fact, $\fF$ is Poisson complete if and only if
$\alpha\wedge\diff\alpha=0$.
\end{example}

\begin{example}\label{exa:poisson_submersion}
Let $\phi:(M,\Pi)\to (B,\pi)$ be any Poisson submersion between two Poisson manifolds. Then the connected
components of the fibers of $\phi$ define a Poisson complete foliation $\fF$ of $M$. This follows from the fact
that the first integrals of $\fF$ are locally of the form $g\circ\phi$, with $g\in C^\infty(B)$, and functions of
this form are closed under the Poisson bracket since $\phi$ is a Poisson map.
\end{example}
\begin{example}
As a particular example of the previous one, consider on $\bbR^2$, with coordinates $(x,y)$, the following Poisson
structure:
$$ \Pi:=x\frac{\partial }{\partial x}\wedge \frac{\partial }{\partial y}\;.$$ The projection on the x-axis,
$(x,y)\mapsto x$ is a Poisson map, so the foliation by vertical lines is Poisson complete. On the other hand, this
foliation is not isotropic since the Poisson tensor vanishes on the vertical line $x=0$, so on points of this line
the inclusion $T\fF\subset \Pi^\sharp(T\fF)^\circ$ does not hold.
\end{example}

\begin{example}\label{exa:coisotropic}
A foliation $\fF$ of a Poisson manifold $(M,\Pi)$ is said to be \textbf{coisotropic} if $\Pi^\sharp(T\fF)^\circ
\subset T\fF$. A necessary and sufficient condition for a foliation $\fF$ of $M$ to be coisotropic is that every
pair of first integrals of the foliation is in involution. Thus, coisotropic foliations are Poisson complete.
\end{example}

%
%

We give in the following proposition a characterization of Poisson complete foliations.
\begin{prop}
\label{prop:Poisson foliations}
Let $\fF$ be an $r$-dimensional foliation of a Poisson manifold $(M,\Pi)$ of dimension $n$. The following
statements are equivalent:
\begin{enumerate}
\item[(i)] $\fF$ is Poisson complete;
\item[(ii)] $(T\fF)^\circ$ is a Lie subalgebroid of $T^*M$.
\end{enumerate}
For any foliation $\fF$ on $(M,\Pi)$ satisfying these conditions, the singular distribution
$\Pi^\sharp(T\fF)^\circ$ is integrable.
\end{prop}
%
%
\begin{proof}
We first recall how the Poisson structure on $M$ makes $T^*M$ into a Lie algebroid (see \cite{Cannas_weinstein} for
background and details). For sections $\a,\b\in\Omega^1(M)$ their Lie bracket is defined by
\begin{equation}\label{eq:bracket_on_forms}
  \lb{\a,\b}:=\L_{\Pi^\sharp(\a)}\b-\L_{\Pi^\sharp(\b)}\a-\diff(\Pi(\a,\b))\;.
\end{equation}%
For (local) sections $f_i\diff g_i$, where $f_i$ is a smooth function, (\ref{eq:bracket_on_forms}) amounts to
\begin{equation}\label{eq:bracket_on_forms_2}
  \lb{f_1\,\diff g_1,f_2\,\diff g_2}=f_1f_2\,\diff\pb{g_1,g_2}+f_1\pb{g_1,f_2}\diff g_2-f_2\pb{g_2,f_1}\diff g_1\;.
\end{equation}%
The anchor of the Lie algebroid $T^*M$ is the map $\Pi^\sharp:T^*M\to TM$.  Let $g_1$ and $g_2$ be two (local)
first integrals of $\fF$ and suppose that $(T\fF)^\circ$ is a Lie subalgebroid of $T^*M$. Then (\ref
{eq:bracket_on_forms_2}) says that $\diff\pb{g_1,g_2}$ is a section of $(T\fF)^\circ$, which means that
$\pb{g_1,g_2}$ is a first integral of $\fF$. This shows that $(ii)$ implies~$(i)$. The converse implication also
follows at once from (\ref {eq:bracket_on_forms_2}) upon using that every section of $(T\fF)^\circ$ is locally of
the form $\sum_i f_i\diff g_i$, where each $g_i$ is a first integral of $\fF$ and the $f_i$ are arbitrary
functions.


The final claim is a consequence of  (ii) because for any Lie algebroid the image of the anchor map is an
integrable (possibly singular) distribution.  
\end{proof}

\begin{prop}
\label{prop:NCIS_bis}
Suppose that $\fF$ is an $r$-dimensional foliation of a Poisson manifold $(M,\Pi)$.
\begin{enumerate}
  \item[(1)] If $(M,\Pi,\fF)$ is an abstract NCI system then $\fF$ is both Poisson complete and isotropic.
  \item[(2)] If $\Pi$ is regular and $\fF$ is both Poisson complete and isotropic, then $(M,\Pi,\fF)$ is an
    abstract NCI system.
\end{enumerate}
\end{prop}

\begin{proof}
(1) Poisson completeness and isotropy of a foliation are local properties, hence they can be proven (easily) by
using Proposition \ref{prp:NCI_local}. Again, we give a direct (easy) proof. Let $m$ be an arbitrary point of $M$
and on a neighborhood $U$ of $m$, let $f_1,\dots,f_r$ be first integrals of $\fF$ whose Hamiltonian vector fields
span $T\fF$. If $g$ and $h$ are first integrals of $\fF$ on $U$, we have in view of the Jacobi identity for $\Pi$:
  $$
    \X_{f_i}(\pb{g,h})=\pb{\X_{f_i}(g),h}+\pb{g,\X_{f_i}(h)}=0\;,
  $$ 
for $i =1,\dots, r$. This shows that the Poisson bracket $\pb{g,h}$ is a local first integral of $\fF$, so $\fF$ is
Poisson complete. Also, since each $f_i$ is a first integral of $\fF$, each $\diff f_i$ is a section of
$(T\fF)^\circ$ and the fact that $T\fF$ is spanned by $\X_{f_1},\dots,X_{f_r}$ implies that $T\fF\subset
\Pi^\sharp(T\fF)^\circ$, so $\fF$ is isotropic.

(2) If $\fF$ is isotropic then $T\fF\subset\Pi^\sharp(T\fF)^\circ$, so that $\Ker \Pi^\sharp\subset
(T\fF)^\circ$. Since $\Pi$ is regular,
\begin{equation*}
  \frac{(T\fF)^\circ}{\Ker \Pi^\sharp}\simeq \Pi^\sharp(T\fF)^\circ
\end{equation*}%
is a (regular) distribution, whose rank is $\Rk(\Pi)-r$. It is generated by the Hamiltonian vector fields
$\Pi^\sharp(\diff f)$ with $f$ a first integral of $\fF$, and these functions are closed under the Poisson bracket,
by Poisson completeness. According to Proposition \ref{prop:Poisson foliations}, this implies that
$\Pi^\sharp(T\fF)^\circ$ is integrable, leading to a foliation $\fG$. If $g$ is a first integral of $\fG$, then
$\X_g$ is tangent to $\fF$. Indeed, if $f$ is a first integral of $\fF$, then $\Pi^\sharp(\diff f)$ is tangent to
$\fG$, so that $\diff f(\Pi^\sharp(\diff g))=-\diff g(\Pi^\sharp(\diff f))=0$.  Note that $\fG$ is contained in the
symplectic foliation of $\Pi$ and has dimension $\Rk(\Pi)+r$. Hence, for any point $m$ of~$M$, we can choose
functions constant on $\fG$ such that at the point $m$:
$$ \diff_m f_1\wedge \dots \wedge \diff_m f_r\ne 0 \quad \hbox{and} \quad \X_{f_1}|_m\wedge \cdots
\wedge\X_{f_r}|_m\ne 0\;. $$
Hence, in a some neighborhood of $m$, the functions $f_1,\dots,f_r$ are constant on $\fF$ and
their Hamiltonian vector fields generate $T\fF$.  This shows that $(M,\Pi,\fF)$ is an abstract NCI system.
\end{proof}

We refer to Section \ref{ex:not_NCI} for an example which shows that an isotropic Poisson complete foliation is not
necessarily an abstract NCI system.

\subsection{Momentum map and Poisson structure on the leaf space}\label{par:NCI_with_mom} 
When the leaf space $B$ of an abstract NCI system $(M,\Pi,\fF)$ is a smooth manifold (i.e., when the holonomy of
$\fF$ is trivial), the leaves of $\fF$ are the connected components of the fibers of the quotient map $\phi: M\to
B$, which is a fibration (with connected fibers). As we will see below (in Proposition \ref{prop:NCIS_ter}), $B$
carries in this case a unique Poisson structure $\pi$ for which $\phi$ is a Poisson map.
\begin{defn}
  We say that an abstract NCI system $(M,\Pi,\fF)$ has a \textbf{momentum map} $\phi:M\to B$ if the leaf space $B$
  of $\fF$ is a (smooth, Hausdorff) manifold. By a small abuse of language, we usually simply speak of an NCI system
  $(M,\Pi)\stackrel{\phi}{\to }(B,\pi)$.

\end{defn}
\begin{prop}
\label{prop:NCIS_ter}
Let $(M,\Pi)\stackrel{\phi}{\to }(B,\pi)$ be an NCI system of dimension~$n$ and rank~$r$. We denote the foliation
on $M$, defined by the fibers of $\phi$, by $\fF$.
\begin{enumerate}
\item[(1)] There exists a unique Poisson structure $\pi$ on $B$ such that $\phi:(M,\Pi)\to (B,\pi)$ is a Poisson
  map;
\item[(2)] Let $f$ be a (local) function, whose Hamiltonian vector field is tangent to the leaves of $\fF$. The
  smooth function $g$ on $B$, defined by $f:=g\circ\phi$, is a (local) Casimir function of $\pi$ (in the
  terminology of Section \ref{par:holonomic}, $g$ is a Cas-basic function);
\item[(3)] For every $m\in M$, $\Rk(\pi_{\phi(m)}) = \Rk(\Pi_m) - 2r $.
\end{enumerate}
\end{prop}
\begin{proof}
Since $\phi$ is a submersion with connected fibers, the smooth functions on~$B$ can be identified with the (global)
first integrals of $\fF$ upon identifying $h\in C^\infty(B)$ with~$h\circ\phi\in C^\infty(M)$. Thus, the Poisson
completeness of $\fF$ leads to (1). It also implies that if $f$ is a function whose Hamiltonian vector field is
tangent to $\fF$, so that $f$ is a first integral of $\fF$, we can write $f$ as $g\circ\phi$ for some function $g$
on $B$. For $h\in C^\infty(B)$ we have that $\pb{h,g}_B\circ\phi=\pb{h\circ\phi,f}=X_f(h\circ\phi)=0$, since
$h\circ\phi$ is a first integral of $\fF$. This shows that $g$ is a Casimir function of~$\PB_B=\pi$, which is the
content of (2).

Let $m\in M$. On a neighborhood of $m$ we can choose functions $f_1,\dots,f_r$ whose Hamiltonian vector fields span
$T\fF$. In view of (2) the functions $g_i$, defined on a neighborhood of $\phi(m)$ by $f_i=g_i\circ\phi$, are
Casimirs of $\pi$.  We denote the differentials of these functions at $m$ and at $\phi(m)$ by $\a_i:=\diff_m f_i$
and $\xi_i:=\diff_{\phi(m)}g_i$. Since the functions $f_i$ are in involution with respect to~$\Pi$, their
(independent) differentials satisfy $\Pi_m(\a_i,\a_j)=0$ for $1\leqs i,j\leqs r$.  They can be completed into a
basis $\a_1,\dots,\a_r$, $\eta_1,\dots,\eta_r$, $\rho_1,\dots,\rho_{n-2r}$ for $T_m^*M$, and since $\Pi_m$ is
skew-symmetric, this can be done such that the matrix of $\Pi_m$ with respect to this basis is
$$
  \left(\begin{array}{ccc}
    0&I_r&0\\
    -I_r&0&0\\
    0&0&Z
  \end{array}
  \right)
$$
where $Z_{ij}=\Pi_m(\rho_i,\rho_j)$. Each one of the $\rho_i$ belongs to $(T_m\fF)^\circ$, since
$\can{\rho_i,\Pi_m^\sharp(\alpha_j)}=\Pi_m(\alpha_j,\rho_i)=0$ for all $j=1,\dots,r$ and since the vectors
$\Pi_m^\sharp(\alpha_j) $ span $T_m\fF$. Therefore, there exist $\s_1,\dots,\s_{n-2r}\in T^*_{\phi(m)}B$ such that
$\rho_i=\phi^*(\s_i)$. In terms of the basis $\xi_1,\dots,\xi_r,\s_1,\dots,\s_{n-2r}$ for $T^*_{\phi(m)}B$, the
matrix of $\pi_{\phi(m)}$ takes the form
$$
  \left(
  \begin{array}{ccc}
    0&0\\
    0&Z
  \end{array}
  \right)
$$
so that the rank of $\pi_{\phi(m)}$ is $2r$ less than the rank of $\Pi_m$, as asserted in~(3).
\end{proof}
\begin{remark}
  Suppose that $(M,\Pi,\F)$ is a regular NCI system of dimension~$n$ and rank $r$ with connected fibers, i.e., the
  fibers of $\F$ are connected. Denoting by $\fF$ the associated foliation and by $B\subset\bbR^{n-r}$ the image of
  $\F$, the abstract NCI system $(M,\Pi,\fF)$ has a momentum map, which is $\F:M\to B$.
\end{remark}

\begin{remark}
Despite the terminology, an abstract NCI system is in general not integrable by quadratures, but an abstract NCI
system with momentum map is. The proof of this fact is essentially the same as in the case of a Liouville
integrable system on a Poisson manifold, see \cite[Sect.\ 4.2]{adlermoerbekevanhaecke2004}. 
\end{remark}

Every abstract NCI system admits a (foliated) atlas, consisting of NCI systems (in the sense of Definition
\ref{def:non-com}), hence it admits locally a momentum map. We will show this in the following proposition. First
we recall (for example from \cite[Ch.~1]{candel}) that an $r$-dimensional foliation $\fF$ of a manifold $M$ of
dimension $n$ can be specified by a \textbf{regular foliated atlas} $(U_\a,\psi_\a\times\phi_\a)_{\a\in I}$: the
$(U_\a)_{\a\in I}$ form an open cover\footnote{The cover can be chosen subordinate to any given open cover of $M$.}
of $M$ and each $\phi_\a$ is a submersion $\phi_\a:U_\a\to\bbR^{n-r}$, whose fibers define the leaves of $\fF$
locally. Moreover, these submersions $\phi_a$ are linked by (unique) diffeomorphisms $\phi_{\a\b}:\phi_\b(U_\a\cap
U_\b)\to \phi_\a(U_\a\cap U_\b)$ such that:
$$ \phi_{\a\b}\circ \phi_\b|_{U_\a\cap U_\b}=\phi_\a|_{U_\a\cap U_\b}\;. $$ 

\begin{prop}
\label{prop:ANCIS}
Let $\fF$ be an $r$-dimensional foliation of a Poisson manifold $(M,\Pi)$ of dimension $n$. The following
statements are equivalent:
\begin{enumerate}
\item[(i)] $\fF$ is an abstract NCI system;
\item[(ii)] $\fF$ admits a regular foliated atlas $(U_\a,\psi_\a\times\phi_\a)$ consisting of NCI systems
  $(U_\a,\Pi\vert_{U_\a},\phi_\a)$ of rank $r$.
\end{enumerate}
\end{prop}

\begin{proof}
The implication (ii) $\Rightarrow$ (i) is straightforward because the foliation defined by a regular NCI system is
an abstract NCI system and because being an abstract NCI system is a local property. Thus, let us suppose
that~$\fF$ is an abstract NCI system on $(M,\Pi)$.
We choose a regular foliated cover $(U_\a)_{\a\in I}$ of~$M$ subordinate to a cover $(U_\b)_{\b\in J}$ having the
property that on each open subset $U_\b$ there exist~$r$ first integrals of $\fF$ (restricted to $U_\b$) whose
Hamiltonian vector fields span~$T\fF$ at every point of~$U_\b$. Let $\a\in I$; we show that $(U_\a,\Pi_{\vert
  U_\a},\psi_\a\times\phi_\a)$ is a regular NCI system. Since the leaves of $\fF$, restricted to $U_\a$, are the
leaves of the foliation of $U_\a$, defined by $\phi_\a$, we may identify local first integrals of~$\fF$, defined on
a neighborhood of a point of $U_\a$ with local first integrals of the foliation defined by the submersion
$\phi_\a$. By construction, there exist first integrals $f_1,\dots,f_r$ on $U_\a$ whose Hamiltonian vector fields
are independent in every point of $U_\a$ (they span $T\fF$ on $U_\a$), in particular their differentials are
independent in every point of $U_\a$. Since $\phi_\a$ is a submersion, there exist extra first integrals
$f_{r+1},\dots,f_s$ of $\fF$, such that $\diff f_1\we\dots\we \diff f_s\neq0$ on~$U_\a$. We have that
$\{f_i,f_j\}=0$ for $1\leqs i\leqs r \hbox{ and } 1\leqs j \leqs s)$, so that
$(U_\a,\Pi\vert_{U_\a},(f_1,\dots,f_s))$ is a regular NCI system, hence also
$(U_\a,\Pi\vert_{U_\a},\psi_\a\times\phi_\a)$.
\end{proof}

\subsection{The semi-local structure of abstract NCI systems with momentum map in the neighborhood of a compact 
fiber}

The existence of local action-angle variables, proved in full generality in \cite{LMV}, can be translated into the
following result, stating that Example \ref{ex:can_model_semi} gives the semi-local model of an abstract NCI system
$(M,\Pi,\fF)$ with a {momentum map}, in the neighborhood of a compact fiber:

\begin{thm}[Semi-local model of an  NCI system with momentum map]\label{thm:aa_semi_local}
Let $(M,\Pi)\stackrel{\phi}{\to}(B,\pi)$ be an NCI system of rank $r=n-s$, where $n$ and $s$ are the dimensions of
$M$ and $B$, respectively, and assume that the fiber $\phi^{-1}(b_0)$ is compact and connected. Then there exist
open neighborhoods $b_0\in U\subset B$ and $0\in V\subset\bbR^s$, a Poisson structure $\pi_0$ on $\bbR^s$ and an
isomorphism $\Psi$ of~NCI systems:
$$
\xymatrix{
(\phi^{-1}(U),\Pi)\ar[rr]^-{\Psi}\ar[d]_-{\phi}& &(\phi_0^{-1}(V),\Pi_0)\ar[d]^-{\phi_0}\\
(U,\pi)\ar[rr]_-{\psi}& &(V,\pi_0)}
$$ 
In this commutative diagram, $\Pi_0,\, \pi_0$ and $\phi_0$ are the Poisson structures and the Poisson map,
defined in Example \ref{ex:can_model_semi}.
\end{thm}

\begin{rem}
\label{rem:polar:foliation}
In the literature (\cite{FA,LMV}) one can find a definition of abstract NCI systems which requires the existence of
a pair of foliations $\fF\subset\fG$ of $(M,\Pi)$ such that $T\fF=\Pi^\sharp(T\fG)^\circ$ (one says that $\fF$ is
polar to $\fG$). For a regular NCI system $(M,\Pi,\F)$ of rank $r$ these foliations are respectively given by the
connected components of the fibers of $\F=(f_1,\dots,f_s)$ and of $\G=(f_1,\dots,f_r)$. The proof of Theorem
\ref{thm:aa_semi_local} given in \cite{LMV} shows that the isomorphism of NCI systems which puts a given NCI system
in a canonical form (providing action-angle coordinates) always respects the foliation~$\fF$, but does not
respect~$\fG$, in general; notice also that although such a foliation~$\fG$ always exists locally, it may not exist
globally (see also Remark~\ref{rem:action:polar:foliation}). For these reasons, we avoid throughout this paper the
assumption of existence of a foliation $\fG$ to which $\fF$ is polar.
\end{rem}


\section{Action variables}\label{sec:action} 

In this section we consider NCI systems with a momentum map, which we write as before as $(M,\Pi)
\stackrel{\phi}{\to }(B,\pi)$. Recall from Section \ref{par:NCI_with_mom} that this means that we have an abstract
NCI system $(M,\Pi,\fF)$, whose leaf space $B$ is a (smooth, Hausdorff) manifold. The latter manifold inherits from
$(M,\Pi)$ a Poisson structure $\pi$ such that the quotient map $\phi:(M,\Pi)\to (B,\pi)$ is a Poisson map. The
foliation $\fF$ is isotropic and Poisson complete. The fibers of $\phi$, which are the leaves of $\fF$, are
connected.

\subsection{The action bundle}

Suppose that we have an NCI system $(M,\Pi) \stackrel{\phi}{\to }(B,\pi)$ of rank~$r$. We construct on $B$ a
canonical vector bundle $E$ of rank~$r$, which is closely related to action variables for the NCI system, as
defined below. To do this, we consider two natural sheaves on $B$ whose quotient essentially represents, pointwise,
the covectors which yield the tangent space to the fibers of $\phi$, upon using the Poisson structure $\Pi$.

Since the bundle map $\pi^\sharp:T^*B\to TB$ may not have constant rank, it is better to view $\pi^\sharp$ as a
sheaf homomorphism $\pi^\sharp\in\HOM(\Omega^1_B,\Vect^1_B)$ from the sheaf $\Omega^1_B$ of differential $1$-forms
on~$B$ to the sheaf $\Vect^1_B$ of vector fields on $B$. Precisely, $\pi^\sharp$ is a homomorphism of sheaves of
$C^\infty_B$-modules: for each open subset $V$ of $B$, we have a $C^\infty_B(V)$-linear map
\begin{equation*}
  \pi^\sharp_V:\Omega^1_B(V)\to \Vect^1_B(V)\;,
\end{equation*}%
which commutes with the restriction maps. The kernel of $\pi^\sharp$ is the subsheaf
$\Ker\pi^\sharp\subset\Omega^1_B$ which to each (non-empty) open subset $V$ of $B$ associates the
$C^\infty_B(V)$-module
$$
  (\Ker\pi^\sharp)(V):=\left\{\omega\in \Omega^1_B(V)\mid \pi_V^\sharp(\omega)=0\right\}\;. 
$$
We also consider another subsheaf $\Ker(\Pi^\sharp\circ\phi^*)\subset \Omega^1_B$ which to each (non-empty) open
subset $V$ of $B$ associates the $C^\infty_B(V)$-module
$$ 
  \Ker(\Pi^\sharp\circ\phi^*)(V):=\left\{\omega\in \Omega^1_B(V)\mid\Pi_{\phi^{-1}(V)}^\sharp(\phi^*\omega)=0\right\}\;. 
$$
Since $\phi$ is a surjective Poisson submersion, $\Ker(\Pi^\sharp \circ\phi^*)(V)\subset(\Ker\pi^\sharp)(V)$, for
every open subset $V$ of $B$.  As a consequence, $\Ker(\Pi^\sharp\circ\phi^*)$ is a subsheaf of $\Ker\pi^\sharp$,
and we can form the quotient sheaf $\cE_B$, which is also a sheaf of $C^\infty_B$-modules on $B$. These sheaves fit
together in the following exact sequence of sheaves on $B$:
$$ 
  \xymatrix{0\ar[r]&\Ker(\Pi^\sharp\circ\phi^*)\ar[r]&\Ker\pi^\sharp\ar[r]& \cE_B\ar[r]&0\;.}
$$ 
Recall from the general theory of sheaves that, for evey open subset $V$ of~$B$, an element of $\cE_B(V)$ is a
collection $(V_i,s_i)_{i\in I}$, where $(V_i)_{i\in I}$ is an open cover of $V$ and $s_i\in(\Ker\pi^\sharp)(V_i)$
for every $i\in I$; these sections are demanded to satisfy $s_i\vert_{V_i\cap V_j}-s_j\vert_{V_i\cap V_j}\in
\Ker(\Pi^\sharp\circ\phi^*)(V_i\cap V_j)$ whenever $V_i\cap V_j\neq\emptyset$. For $\omega\in(\Ker\pi^\sharp)(V)$
we denote its image in $\cE_B(V)$ by $[\omega]$.

The above construction works for any surjective Poisson submersion $\phi:(M,\Pi)\to (B,\pi)$. We show in the
following proposition that for an NCI system of rank $r$ the sheaf $\cE_B$ is the sheaf of sections of a vector
bundle $E\to B$ of rank $r$.
\begin{prop}\label{prp:ranks}
  Let $(M,\Pi) \stackrel{\phi}{\to }(B,\pi) $ be an NCI system of rank $r$. The quotient sheaf $\cE_B:=
  \Ker\pi^\sharp/\Ker(\Pi^\sharp\circ\phi^*)$ is the sheaf of sections of a vector bundle $E$ on $B$ of rank $r$.
  We call $\cE_B$ the {\bf action sheaf} and $E\to B$ the {\bf action bundle} of the NCI system.
\end{prop}
\begin{proof}
We need to show that $\cE_B$ is a locally free sheaf of $C^\infty_B$-modules of rank~$r$. Let $b\in B$ and denote,
as before, by $\fF$ the foliation of $M$ defined by the fibers of $\phi$. According to the definition of an NCI
system and Proposition~\ref{prop:NCIS_ter}~(2) there exist, on a neighborhood $V$ of $b$, Casimir functions
$g_1,\dots,g_r$ of $\pi$ such that $T\fF$ is spanned at each point of $\phi^{-1}(V)$ by $\X_{f_1},\dots,\X_{f_r}$,
where $f_i:=g_i\circ\phi$, for $i=1,\dots,r$. Let $s\in\cE_B(V)$. By definition, $s$ is given by a collection
$(V_i,s_i)_{i\in I}$, where $(V_i)_{i\in I}$ is an open cover of $V$ and
$s_i\in(\Ker\pi^\sharp)(V_i)\subset\Omega^1_B(V_i)$ for every $i\in I$; also $s_i\vert_{V_i\cap
  V_j}-s_j\vert_{V_i\cap V_j}\in \Ker(\Pi^\sharp\circ\phi^*)(V_i\cap V_j)$ whenever $V_i\cap
V_j\neq\emptyset$. Since the vector fields $\Pi^\sharp(\phi^*s_i)$ are tangent to the fibers of $\phi$, there exist
unique smooth functions $\l_{il}$ on $\phi^{-1}(V)$, such that
\begin{equation}\label{eq:pullb}
  \Pi^\sharp(\phi^*s_i)=\sum_{l=1}^r\l_{il}\Pi^\sharp(\diff f_l)\;.
\end{equation}%
We show that these functions are $\phi$-basic (i.e., constant on the fibers of $\phi$). To do this, we show that
$\X_{f_k}(\l_{il})=0$ for $i\in I$ and $k,l=1,\dots,r$. Since $X_{f_k}$ is tangent to $\fF$,
\begin{equation*}
  \lb{\X_{f_k},\Pi^\sharp(\phi^* s_i)}=\Pi^\sharp(\cL_{X_{f_k}}\phi^*s_i)=0\;,  
\end{equation*}%
so that
\begin{equation*}
  \sum_{l=1}^rX_{f_k}(\l_{il})\Pi^\sharp(\diff f_l)=0\;.
\end{equation*}%
This shows our claim because the vector fields $\Pi^\sharp(\diff f_1),\dots,\Pi^\sharp(\diff f_r)$ are linearly
independent at every point of $\phi^{-1}(V)$. Since the fibers of $\phi$ are connected, it follows that there exist
(unique) smooth functions $\s_{il}$ on $V$ such that $\l_{il}=\s_{il}\circ\phi$. Substituted in (\ref{eq:pullb}),
we find that
\begin{equation*}
  \Pi^\sharp\phi^*\(s_i-\sum_{l=1}^r\s_{il}\diff g_l\right)=0\;,
\end{equation*}%
so that $s_i-\sum_{l=1}^r\s_{il}\diff g_l\in\Ker(\Pi^\sharp\circ\phi^*)(V_i)$. For $i,j$ such that $V_i\cap
V_j\neq\emptyset$ we have that $s_i\vert_{V_i\cap V_j}-s_j\vert_{V_i\cap V_j}\in
\Ker(\Pi^\sharp\circ\phi^*)(V_i\cap V_j)$, so that $\s_{il}=\s_{jl}$ on $V_i\cap V_j$ for all $l$. Thus, the
functions $(\s_{il})_{i\in I}$ glue together to a global function $\s_l\in C^\infty_B(V)$ and we can write
$s=\sum_{l=1}^r\s_l[\diff g_l]$ for some unique smooth functions $\s_l$ on $V$, as required.
\end{proof}
For $b\in B$, the fiber $E_b$ of the vector bundle $E\to B$ corresponding to $\cE_B$ can be recovered from $\cE_B$
as
\begin{equation}\label{eq:E_recover}
  E_b=\frac{\cE_B(V)}{C^\infty_b(V)\cE_B(V)}\;,
\end{equation}%
where $C^\infty_b(V)$ stands for the ideal of $C^\infty_B(V)$ containing all smooth functions on $V$ which vanish
at $b$ and $V$ is any open subset of $B$ containing $b$ and such that $\cE_B(V)$ is a free $C^\infty_B(V)$-module.
Let $m$ be any point in the fiber of $\phi$ over $b$. We show that the following sequence of vector spaces is
exact:
\begin{equation}\label{eq:action_exact}
  \xymatrix{0\ar[r]&C^\infty_b(V)\cE_B(V)\ar[r]&\cE_B(V)\ar[r]^-{\rho_b}&
  \displaystyle\frac{\Ker(\pi_{b}^\sharp)}{\Ker(\Pi_m^\sharp\circ\phi^*)}\ar[r]&0\;.}
\end{equation}
To do this, we first show that if $m,m'\in\phi^{-1}(b)$ then $\Ker(\Pi_m^\sharp \circ\phi^*)=\Ker( \Pi_{m'}^\sharp
\circ\phi^*),$ so that the latter space is independent of the choice of $m$ in $\phi^{-1}(b)$. Since the fibers of
$\phi$ are connected it is enough to prove the equality for $m'$ in a neighborhood of~$m$. There exist, in a
neighborhood of $\phi(m)$ in $B$, Casimir functions $g_1,\dots,g_r$ such that the Hamiltonian vector fields of
$f_1:=g_1\circ\phi,\dots,f_r:=g_r\circ\phi$ span $T\fF$ in a neighborhood of $m$ in~$M$. The (local) flows of these
vector fields commute, since $\lb{X_{f_i},X_{f_j}}=-X_{\pb{g_i,g_j}\circ\phi}=0$. These flows therefore define a
(local) action of $\bbR^r$, by Poisson diffeomorphisms, which is transitive in a neighborhood of $m$. In
particular, we obtain a local Poisson diffeomorphism $\Psi$ such that $\phi\circ\Psi=\phi$ and
$\Psi(m)=m'$. It follows that:
$$
  \Pi^\sharp_{m'}\circ \phi^*=\diff_m\Psi\circ \Pi^\sharp_{m}\circ (\diff_{m}\Psi)^*\circ \phi^*=\diff_m\Psi\circ \Pi^\sharp_{m}\circ \phi^*\;.
$$
This implies our claim since $\diff_m\Psi$ is an isomorphism. We can now prove that (\ref{eq:action_exact}) is a
short exact sequence. Since the injectivity of the first arrow and the surjectivity of the last arrow are clear, we
only prove the exactness at $\cE_B(V)$. Let $s$ be an element of~$\cE_B(V)$. As we have seen in the proof of
Proposition~\ref{prp:ranks}, $s$ can be written as $s=\sum_{l=1}^r\s_l[\diff g_l]$ for some unique smooth functions
$\s_l$ on $V$. Exactness then follows from the fact that $\rho_b(s)=\sum_{l=1}^r\sigma_l(b)[\diff_b g_l]$ where, by
a slight abuse of notation, $[\diff_b g_l]$ stands for the class of $\diff_b g_l$ in
${\Ker(\pi_{b}^\sharp)}/{\Ker(\Pi_m^\sharp\circ\phi^*)}$.

The exactness of (\ref{eq:action_exact}), combined with (\ref{eq:E_recover}), provides a natural identification of
$E_b$ with $\frac{\Ker(\pi_{b}^\sharp)}{\Ker(\Pi_m^\sharp\circ\phi^*)}$. As we show next, the Poisson structure
$\Pi$ also induces a natural identification of $E_b$ with $T_m\fF$, which is the tangent space to $\phi^{-1}(b)$ at
$m$, where $m$ is an arbitrary point in $\phi^{-1}(b)$. Indeed, every $\alpha \in E_b $ defines a smooth vector
field $X_\alpha$ on the fiber $\phi^{-1}(b) $ over $b$, by
$$
  X_\alpha(m) := \Pi_m^\sharp(\phi^*\a)
$$ 
for all $m \in \phi^{-1}(b) $. We call $X_\a$ the {\bf action vector field} associated to~$\a$. 

\begin{lemma}\label{lma:fund_vfs}
  Let $(M,\Pi) \stackrel{\phi}{\to }(B,\pi) $ be an NCI system of rank $r$. Let $m\in M$ and denote $b:=\phi(m)\in
  B$.
  \begin{enumerate}
    \item For every $\alpha,\alpha' \in E_b$, the action vector fields $X_\alpha$ and $X_{\alpha'}$ commute;
    \item For every basis $\alpha_1, \dots,\alpha_r $ of $E_b $ the vector fields $X_{\alpha_1},\dots,
      X_{\alpha_r}$ form a basis of $T_m\fF$. In particular, $X_\a$ is nowhere vanishing when $\a\neq0$.
  \end{enumerate}
\end{lemma}
\begin{proof}
On a neighborhood $V$ of $b$ in $B$ there exist Casimirs $g_1,\dots,g_r$ such that their associated vector fields
$\Pi^\sharp(\diff(g_i\circ\phi))$ generate the tangent space to $\fF$ on $\phi^{-1}(V)$. It follows that
$[\diff_b{g_1}],\dots,[\diff_bg_r]$ are independent, hence form a basis for $E_b $ and (2) follows. The vector
fields $X_{g_1\circ\phi},\dots,X_{g_r\circ\phi}$ are tangent to the fibers of $\phi$ over $V$ and they commute, as
we have seen above. In particular, the vector fields $X_\a$ commute.
\end{proof}
In view of item (2) above, the map $E_{\phi(m)}\to T_m\fF$ which sends $\a\in E_{\phi(m)}$ to $X_\a$ is an
isomorphism and we may think of $E_{\phi(m)}$ as being the tangent space to the fiber of $\phi$ at $m$.

The notation $X_\a$ which we introduced above for elements $\a$ of $E_b$ will also used for (local) sections of
$E\to B$: for a section $e\in\cE_B(V)$, the {\bf action vector field} $X_e$ is a vector field which is defined on
$\phi^{-1}(V)$ and it is tangent to the fibers of $\phi$: for $b\in V$, the restriction of $X_e$ to $\phi^{-1}(b)$
is~$X_{e(b)}$. For arbitrary sections $e,e'\in\cE_B(V)$ the vector fields $X_e$ and $X_{e'}$ commute, in view of
Lemma \ref{lma:fund_vfs} (1).

\subsection{Holonomic sections of the action bundle}\label{par:holonomic}
Suppose that we have an NCI system $(M,\Pi) \stackrel{\phi}{\to }(B,\pi)$ of rank~$r$. We denote its action sheaf
by~$\cE_B$. For $V\subset B$, we call an element $e\in\cE_B(V)$ \textbf{locally holonomic} if for every point $b\in
V$, there exists a Casimir function $g$ of $\pi$, defined on a neighborhood $W\subset V$ of $b$ in $B$, such that
$e_{|_W}=[\diff g]$. Notice that when two such neighborhoods $W_1$ and $W_2$ intersect, the Casimir functions $g_1$
and~$g_2$ which define $e$ satisfy $[\diff(g_1-g_2)]=0$ on $W_1\cap W_2$, so that $\phi^*(g_1-g_2)$ is a Casimir
function of $\Pi$ (on $\phi^{-1}(W_1\cap W_2)$). Therefore, we introduce three more sheaves on $B$, by letting for
every open subset $V$ of $B$:
\begin{eqnarray*}
  \Cas_B(V)&:=&\left\{F\in C^\infty(V)\mid F\hbox { is a Casimir function of } \pi\vert_V\right\}\;,\\
  \CASM(V)&:=&\left\{F\in C^\infty(V)\mid F\circ\phi \hbox { is a Casimir function of } \Pi\vert_{\phi^{-1}(V)}\right\}\;,\\
  \cE_B^{0}(V)&:=&\left\{e\in\cE_B(V)\mid e\hbox{ is locally holonomic}\right\}\;.
\end{eqnarray*}
$\Cas_B$ is the sheaf of Casimir functions on $B$, while $\CASM$ is the sheaf of \textbf{Cas-basic functions}, that
is local functions on $B$ whose pullback to $M$ are Casimir functions. Notice that, contrary to the sheaves which
were introduced in the previous subsection, they are simply sheaves of $\bbR$-vector spaces and not of
$C^\infty_B$-modules.  Since $\phi$ is a surjective Poisson morphism, $\CASM$ is included in $\Cas_B$, and the
above argument shows that $\cE_B^0$ is the quotient sheaf $\Cas_B/\CASM$, i.e., the following sequence of sheaves
of vector spaces on $B$ is exact:
\begin{equation}\label{eq:E0_seq}
  \xymatrix{0\ar[r]&\CASM\ar[r]&\Cas_B\ar[r]^-{[\diff\cdot]}& \cE^0_B\ar[r]&0\;.}
\end{equation}
For future use, we give the following exact sequence of sheaves, which derives from the previous one:
\begin{equation}\label{eq:E0_seq_bis}
  \xymatrix{0\ar[r]&\CASM/\bbR\ar[r]&\Cas_B/\bbR\ar[r]^-{[\diff\cdot]}& \cE^0_B\ar[r]&0\;;}
\end{equation}
here, and in all further sheaf contexts, $\bbR$ stands for the sheaf of locally constant functions on the manifold
under consideration, in this case $B$.

For an open subset $V$ of $B$, an element $e$ of $\cE_B(V)$ is called a \textbf{globally holonomic} section if
$e=[\diff g]$ for some Casimir function $g$ on $V$.  We will be particularly interested in globally holonomic
sections which are defined on all of $B$. In order to characterize these sections, we consider the long cohomology
sequence associated to (\ref{eq:E0_seq}), which is given in part by
$$
  \xymatrix{\cdots\ar[r]& H^0(B,\Cas_B)\ar[r]& H^0(B,\cE^0_B)\ar[r]^-{\Obs}& H^1(B,\CASM)\ar[r]&\cdots}
$$
The connecting homomorphism defines a map, which we denote by $\Obs$ and which we call the \textbf{holonomy
  obstruction} (of the NCI system). The locally holonomic elements of $\cE_B(B)$ are precisely the elements of
$H^0(B,\cE^0_B)$, while the globally holonomic elements of $\cE_B(B)$ are the elements in the image of
$H^0(B,\Cas_B)\rightarrow H^0(B,\cE^0_B)$. Exactness of the above long exact sequence leads to the following
proposition.
\begin{prop}\label{prp:holonomic}
  Let $e$ be a global section of $\cE_B$ which is locally holonomic. Then $e$ is globally holonomic if and only
  if $\Obs(e)=0$.
\end{prop}

\subsection{The action lattice bundle and the integral affine structure on the fiber}\label{par:lattice}
We say that an NCI system $(M,\Pi) \stackrel{\phi}{\to }(B,\pi) $ has {\bf compact fibers} when all the fibers of
$\phi$ are compact.  In this case, all vector fields $X_{g\circ\phi}$, with~$g$ a Casimir function of $\pi$,
defined on an open subset of $B$, are complete. In particular, the action vector fields~$X_\a$, with $\a\in E_b$
are complete and we can consider their time~1 flow. In view of Lemma \ref{lma:fund_vfs} the action vector fields
associated to two elements of $E_b$ (with $b\in B$) commute, hence the time~1 flow defines an action of $E_b$ on
$\phi^{-1}(b)$. By the same lemma, the action is locally free, hence transitive (recall that by definition the
fibers of $\phi$ are connected).  It follows that there is for each $b\in E_b$ a canonically defined lattice ${L_b}
\subset E_b $, namely the lattice of all points $\alpha \in E_b$ such that the time~1 flow of~$ X_\alpha$ is the
identity map.  Said differently,~$L_b$ is the subset of all the elements $\a$ of~$E_b$ such that for one
(equivalently, for all) $m\in \phi^{-1}(b)$ the time~1 flow of the action vector field~$X_a$ fixes $m$. We call
$L_b\subset E_b$ the {\bf action lattice} at $b$. As $b$ runs through~$B$, these lattices $L_b$ fit nicely together
in a group bundle $L$ over $B$, with fiber $\bbZ^r$; for the proof of this fact, we refer to
\cite[Sect.\ 3.4]{LMV}.  We call $L$ the {\bf action lattice bundle} of the NCI system.

We will find it convenient to view the local sections of $L\to B$ as a sheaf on $B$, which we denote by $\cL_B$ and
which we call the {\bf action lattice sheaf}. Thus, for any open subset $V$ of $B$ we denote by $\cL_B(V)$ the
space of sections of $L\to B$ over $V$. It is clear that $\cL_B$ is a sheaf of $\bbZ$-modules on $B$:
locally,~$\cL_B$ is isomorphic to the constant sheaf $\bbZ^r$ on $B$. An isomorphism between the restrictions of
$\cL_B$ and $\bbZ^r$ to $V\subset B$ is called a {\bf trivialization} of $\cL_B$ on~$V$. Such as isomorphism is
defined by $r$ sections of $\cL_B$ over $V$.

We can now define the notion of action variables in terms of the above terminology. Let $V$ be an open subset of
$B$. We say that an $r$-tuple $(p_1, \dots, p_r)$ of functions on $V$ are a set of {\bf local action variables} (on
$V$) if $ [\diff p_1], \dots, [\diff p_r]$ define a trivialization of ${\cL}_B$ on~$V$. In view of
Proposition~\ref{prop:NCIS_ter} (2), local action variables are (local) Casimir functions of $(B,\pi)$. Local
action variables on $V=B$ are called {\bf global action variables}. Since, as we pointed out above, we can identify
functions on $B$ with functions on~$M$ which are constant on the fibers of $\phi$, we will also call the functions
$\phi^*p_i$ a set of (local or global) action variables.
	
\begin{rem}
  By construction, if $(p_1, \dots, p_r)$ are a set of action variables on~$V$, then the Hamiltonian vector fields
  of $\phi^* p_1, \dots, \phi^* p_r$ are periodic of period one and they commute; in particular $\phi^* p_1, \dots,
  \phi^*p_r$ are the components of a momentum map of a $\bbT^r$ action on $V$. These properties justify the
  terminology \emph{action variables}.
\end{rem}	
\noindent Theorem~\ref{thm:aa_semi_local} implies the following results.
\begin{prop}\label{prp:action_vars_exist}
  Let $(M,\Pi) \stackrel{\phi}{\to} (B ,\pi) $ be an NCI system with compact fibers. 
  \begin{enumerate}
    \item[(1)] Local action variables exist on a neighborhood $V$ of every point $b\in B$;
    \item[(2)] $\cL_B$ is a subsheaf of $\cE^0_B$, where both sheaves are viewed as sheaves of $\bbZ$-modules. Said
      differently, if $V$ is an open subset of $B$ and $\ell\in\cL_B(V)$, then $\ell$ is locally holonomic.
\end{enumerate}
\end{prop}
The following theorem gives a cohomological condition for the existence of global action variables.
\begin{thm}\label{thm:action_var_obs}
  Let $(M,\Pi) \stackrel{\phi}{\to }(B,\pi) $ be an NCI system of rank $r$ with compact fibers. The following
  properties are equivalent:
    \begin{enumerate}
      \item[(i)] There exists a set of  global action variables;
      \item[(ii)] The action lattice sheaf ${\cL_B} $ admits a (global) trivialization and every global section
        $\ell$ of ${\cL_B} $ satisfies $\Obs(\ell)=0$.
    \end{enumerate}
\end{thm}
\begin{proof}
Let $(p_1,\dots,p_r)$ be a set of global action variables for the NCI system.  By definition, $[\diff p_1], \dots,
[\diff p_r]$ define a trivialization of $\cL_B$ on $B$, hence every global section $ \ell$ of $\cL_B$ is of the
form
$$
  \ell = \sum_{i=1}^r n_i [\diff p_i] 
$$ 
for some integers $n_1,\dots,n_r$. This implies that $\ell =[\diff(\sum_{i=1}^r n_i p_i)]$ is in the image of
$\Cas_B(B)\to\cE^0_B(B)$, so that $\Obs(\ell)=0$.  This proves that (i) implies~(ii).
        
Conversely, suppose that $\ell_1,\dots,\ell_r$ define a trivialization of $\cL_B$ on $B$ and that $\Obs(\ell_i)=0$
for $i=1,\dots,r$. According to Proposition \ref{prp:holonomic}, there exist Casimir functions $p_1,\dots,p_r$ such
that $\ell_i = [\diff p_i]$ for $i=1,\dots,r$.  By definition, $(p_1, \dots,p_r)$ is a set of global action
variables for the NCI system.
\end{proof}

\begin{rem}\label{rem:LB_trivial}
  Notice that saying that ${\cL_B} $ admits a (global) trivialization is equivalent to saying that $ M \to B$ is a
  principal $\bbT^r$-bundle; it is also equivalent to saying that the class defined by $\cL_B$ in
  $H^1(M,GL_r(\bbZ))$ is trivial.
\end{rem}

\subsection{Action foliations}
\label{sec:action:foliation}
%
Let $(M,\Pi) \stackrel{\phi}{\to }(B,\pi) $ be an NCI system of rank $r$ with compact fibers. A foliation $\fA$ of
$B$ is said to be an {\bf action foliation} of the NCI system when ${\fA}$ is defined in the neighborhood of every
point by local action variables. It means that on a neighborhood $V$ of any point $b\in B$ we can find functions
$p_1,\dots,p_r$ such that 
\begin{enumerate}
  \item[(1)] $[\diff p_1],\dots,[\diff p_r]$ define a trivialization of $\cL_B$ on $V$;
  \item[(2)] The foliation $\fA$, restricted to $V$, is defined by $p_1,\dots,p_r$.
\end{enumerate}
Obviously, the foliation defined by action variables is an action foliation, but the converse is false in general,
as we will see. 

\begin{remark}
\label{rem:action:polar:foliation}
The foliation $\fF$ of the NCI system is polar to the pullback $\phi^{-1}(\fA)$ of any action foliation $\fA$ (see
Remark \ref{rem:polar:foliation}). Note, that if $\fF$ is polar to some foliation $\fG$, then $\fG$ is the pullback
of a foliation $\phi(\fG)$, but this foliation, in general, will fail to be an action foliation. One can show that
this is the case if and only if $\fG$ is locally given around its leaves by the kernel of basic closed 1-forms
$\alpha_1,\dots,\alpha_r\in\Omega^1(M)$ with the property that the the vector fields
$\Pi^\sharp(\alpha_1),\dots,\Pi^\sharp(\alpha_1)$ have all their orbits periodic with period 1. Hence, the
existence of an action foliation requires the existence of a polar foliation of a very special nature.
\end{remark}

We denote by $\Cas_\fA$ the sheaf of local first integrals of $\fA$; the notation is motivated by
the first item in the following proposition:


\begin{prop}\label{prop:actionfoliation}
  Let $(M,\Pi) \stackrel{\phi}{\to }(B,\pi) $ be an NCI system of rank $r$ with compact fibers. Suppose that it has
  an action foliation ${\fA}$. Then the following properties are satisfied:
  \begin{enumerate}
    \item[(1)] $\Cas_\fA$ is a subsheaf of $\Cas_B$; said differently, $\fA$ contains the symplectic foliation of
        $(B,\pi)$;
%
    \item[(2)] ${\fA} $ is a transversely integral affine foliation.
  \end{enumerate}
\end{prop}
\begin{proof}
Item (1) follows from the fact that $\fA$ is locally defined by action variables, which are (local) Casimir
functions of $\pi$. 
%
In order to prove (2), consider a cover of $B$ by open sets on which $\fA$ is defined by local action
variables. Let $V$ and $V'$ be two intersecting subsets of the cover and let $(p_1, \dots,p_r)$
(resp.\ $(p_1',\dots,p_r')$) be a set of action variables on $V$ (resp.\ on $V'$) which define $\fA$. Then we can
write on a connected neighborhood $W$ of any $b\in V\cap V'$ the functions $p'_1, \dots,p'_r$ in terms of
$p_1,\dots,p_r$. Taking the differential, we get
$$ 
  \diff p_i'=\sum_{k=1}^r\pp{p_i'}{p_k}\diff p_k\;,\qquad (i=1,\dots,r)\;.
$$
Since both $([\diff p_1],\dots,[\diff p_r])$ and $([\diff p_1'],\dots,[\diff p_r'])$ define a trivialization
of~${\cL}_B$ on $W$, the above relations imply that the functions $a_{ij}:=\pp{p_i'}{p_j}$ are constant and take
values in $\bbZ$, for all $i,j =1, \dots, r$. Since $W$ is connected, it follows that each one of the functions
$p_1, \dots,p_r$ is, up to real a constant, a linear combination with integral coefficients of the functions $p_1,
\dots,p_r$; this is precisely the property which defines transversely integral affine foliations.
\end{proof}

\begin{remark}\label{rem:action_foliation_symplectic}
  When $M$ is symplectic, Proposition \ref{prop:NCIS_ter} (2) implies that the Poisson structure $\pi$ on $B$ is
  regular, with symplectic leaves of dimension $\dim B-r$. Every set of local action variables defines the
  symplectic foliation, hence there exists precisely one action foliation, which coincides with the symplectic
  foliation. 
\end{remark}

We will only analyse the obstruction to the existence of an action foliation when $\cL_B$ admits a trivialization
over $B$, Associated to the following short exact sequence of sheaves on $B$:
$$  
  0 \to \bbR \to \CASM \to  \CASM/\bbR  \to 0 \;, 
$$
there is the long exact sequence
\begin{equation} \label{eq:longexactDiviserR}  
  \cdots \to H^1(B, \bbR ) \to  H^1(B,\CASM) \to  H^1(B,\CASM/\bbR)  \to \cdots 
\end{equation}
We say that a class in $H^1(B,\CASM)$ is {\bf representable by constants} if it lies in the image of $H^1(B, \bbR )
\to H^1(B,\CASM)$, or, equivalently, in the kernel of $ H^1 (B,\CASM) \to H^1 (B,{\CASM/\bbR}) $.  A class is
representable by constants if and only if it can be represented by a cocycle valued in locally constant functions,
hence the name.

\begin{prop}\label{cor:actionfoliation}
  Let $(M,\Pi) \stackrel{\phi}{\to }(B,\pi) $ be an NCI system of rank $r$ with compact fibers. Suppose that its
  action lattice sheaf $\cL_B$ admits a trivialization on $B$, defined by sections $\ell_1,\dots,\ell_r$ of $\cL_B$
  over $B$. Then the following conditions are equivalent:
  \begin{enumerate}
    \item[(i)] There exists a global action foliation for the NCI system;
    \item[(ii)] For $i=1,\dots,r$, the class $\Obs (\ell_i) \in H^1(B,\CASM)$ is representable by constants.
  \end{enumerate}
\end{prop}	
\begin{proof}
The connecting morphism $\Obs$ of the exact sequence (\ref{eq:E0_seq}) and the connecting morphism $\delta$ of the
exact sequence (\ref{eq:E0_seq_bis}) are related through the following commutative diagram:
$$ 
 \xymatrix{&H^0(B,\cE_B^0)\ar[d]^-{\Obs}\ar[dr]^-{\delta}&\\  
           H^1(B, \bbR ) \ar[r] &  H^1(B,\CASM) \ar[r] &  H^1(B,\CASM/\bbR)}
$$ 
Since the horizontal line of this diagram is exact, $\Obs(\ell_i)$ is representable by constants if and only if
$\delta(\ell_i)=0$. 
%

Suppose that there exists a global action foliation $\fA$. Then there in the neighborhood of every point of $B$
Casimir functions $p_1,\dots,p_r$, such that $[\diff p_i]=\ell_i$ for $i=1,\dots,r$. As in the proof of Proposition
\ref{prop:actionfoliation} the functions $p_i$ and $p_i'$ differ on overlapping opens only by locally constant
functions, hence the cocycle which is defined by $\ell_i$ is trivial in $H^1(B,\CASM)$, i.e., $\delta(\ell_i)=0$,
so that $\Obs(\ell_i)$ is representable by constants. This shows that (i) implies~(ii). 

Suppose now that each $\Obs(\ell_i)$ is representable by constants. Then there exists a cover of $B$ by open
subsets $(U_j)_{j\in J}$ and Casimir functions $p_{1j},\dots,p_{rj}$ on each $U_j$, such that for every $i=1,\dots,r$,
\begin{enumerate}
  \item[(1)] $[\diff p_{ij}]=\ell_i$ on $U_j$, for all $j\in J$;
  \item[(2)] On non-empty overlaps $U_j\cap U_k$, which are supposed connected, $p_{ij}-p_{ik}$ is constant.
\end{enumerate}
The first condition implies that for fixed $j\in J$ the functions $p_{1j},\dots,p_{rj}$ define an action foliation
on $U_j$, while the second condition implies that the action foliations on $U_j$ and $U_k$ concide on $U_j\cap
U_k$, hence define a global action foliation on $B$. This shows that (ii) implies (i).
\end{proof}

\begin{remark}
Proposition \ref{cor:actionfoliation} can be generalized to the case where the lattice sheaf is not trivial. This
can be done as follows: let $(M,\Pi) \stackrel{\phi}{\to }(B,\pi) $ be an NCI system of rank $r$ with compact
fibers and let $\cL_B$ denote its lattice sheaf. The following statements are then equivalent:
  \begin{enumerate}
    \item[(i)] There exists an  action foliation for $(M,\Pi) \stackrel{\phi}{\to }(B,\pi) $;
    \item[(ii)] The cohomology class $\Obs ({\cL_B}) \in  H^1(B,\HOM_\bbZ(\cL_B,\CASM/\bbR))$ vanishes.
  \end{enumerate}
Let us define the class and the cohomology space that appear in {(ii)}. For any sheaf of abelian groups ${\mathcal
  F}$ over $B$, we denote by $\HOM_\bbZ(\cL_B,{\mathcal F})$ the sheaf whose sections over an open subset $U
\subset B$ is the set of all group morphisms from $\cL_B(U)$ to ${\cF}(U)$; thus, $\HOM_\bbZ(\cL_B,{\mathcal F})$
is itself a sheaf of abelian groups. Applying to the exact sequence~(\ref{eq:E0_seq_bis}) the exact functor
$\HOM_\bbZ(\cL_B, \cdot)$ yields an exact sequence:
\begin{equation}\label{eq:exactsequenceHom}
  0 \to \HOM_\bbZ(\cL_B,\CASM/\bbR) \stackrel{}{\to} \HOM_\bbZ(\cL_B,\Cas_B/\bbR) \stackrel{[\diff \cdot]}{\to} \HOM_\bbZ(\cL_B,\cE_B^0) \to 0 \;.
\end{equation} 
Now, the canonical inclusion $\cL_B \hookrightarrow \cE_B^0$ can be seen as an element in $ H^0(B,\HOM_\bbZ
(\cL_B,\cE_B^0)) $, to which the connecting morphism of (\ref{eq:exactsequenceHom}) can be applied, giving a class
in $ H^1(B,\HOM_\bbZ(\cL_B,\CASM/\bbR))$, which we denote by $\Obs ({\cL_B})$. 

The proof of the equivalence between (i) and (ii) follows essentially the same lines as the proof of proposition
\ref{cor:actionfoliation}, upon noticing that $\delta(\imath_\cL)=0$ is tantanamount to the existence of a sheaf
homomorphism $\jmath_{\cL} $ from $\cL_B $ to $\Cas_B/\bbR$ which makes the following diagram commutative:
\begin{equation}\label{eq:commutativediagramm=actionfoliation}
\xymatrix{\cL_B\ar@{^{(}-->}[rr]^-{\jmath_\cL}\ar@{^{(}->}[dr]_-{\imath_\cL}&&\Cas_B/\bbR\ar@{_{(}->}[ld]^-{[\diff\cdot]}\\
          &\cE_B^0&}
\end{equation}
while the existence of $\jmath_{\cL} $ can be checked to be equivalent to the existence of an action foliation.
\end{remark}

%
\section{Angle variables and transverse structure}\label{sec:angle} 

In this section, we suppose that we have an NCI system $(M,\Pi) \stackrel{\phi}{\to }(B,\pi)$ with compact
fibers. As before, we denote its action lattice sheaf by $\cL_B$.

\subsection{Angle variables}\label{par:angle_variables}

We first define the notion of angle variables.

\begin{defn}\label{def:angles_variables}
  Let $(e_1, \dots, e_r )$ be a trivialization of $\cL_B(V)$ where $V$ is some open subset of $B$. An $r$-tuple of
  $\bbR/\bbZ$-valued functions $(\theta_1,\dots,\theta_r)$ defined on $\phi^{-1}(V)$ is called a set of {\bf local
    angle variables} on $\phi^{-1}(V)$, adapted to $(e_1, \dots, e_r )$, if
  \begin{equation}\label{eqn:angle_def}
    \pb{\theta_i,\theta_j}=0\;,\qquad X_{e_i}(\theta_j)=\delta_{i,j}\;,
  \end{equation}
  for all $1\leqs i,j\leqslant r$. 
\end{defn}

Notice that, given a set of local angle variables, the trivialization of $\cL_B(V)$ with respect to which it is
adapted is uniquely determined by it, so we may speak of local angle variables without specifying a (local)
trivialization of~$\cL_B$. As a consequence, given an $r$-tuple of $\bbR/\bbZ$-valued functions
$(\theta_1,\dots,\theta_r)$ on $M$, which are local angle variables in the neighborhood of every point of $B$,
there exists a (global) trivialization $(e_1,\dots,e_r)$ of $\cL_B(B)$ such that $(\theta_1,\dots,\theta_r)$ are
angle variables on $M$, adapted to it. We then call $(\theta_1,\dots,\theta_r)$ {\bf global angle variables}.

The following proposition is a corollary of the local action-angle theorem (Theorem~\ref{thm:aa_semi_local}):
\begin{prop}\label{prp:angels_exist}
Every point $b \in B$ is contained in an open neighborhood~$V$ such that there exists a trivialization
$e=(e_1,\dots, e_r)$ of $\cL_B(V)$ and a set of local angle variables on $\phi^{-1}(V) $ adapted to $e$.
\end{prop}

In order to show how two different sets of local angle variables are related, we first construct $r$ vector fields
$Y_{\theta_i}$ on $V\subset B$ which represent the Hamiltonian vector fields, associated to a set of local angle
variables\footnote{We will see in Section \ref{par:transverse_foliation} that these vector fields define an
integrable distribution of rank $r$ which depends only on the foliation, defined by the angle variables and which
is transverse to every local action foliation.}  of the lattice sheaf.
\begin{prop}\label{prop:angles_vectorfield_quotient}
  Let $V$ be an open subset of $B$ and suppose that $(\theta_1, \dots, \theta_r )$ is a set of local angle
  variables on $\phi^{-1}(V)$, adapted to some trivialization $(e_1,\dots,e_r) $ of $\cL_B(V)$. The Hamiltonian
  vector fields $X_{\theta_1}, \dots,X_{\theta_r}$ are $\phi$-related to commuting Poisson vector fields
  $Y_{\theta_1}, \dots,Y_{\theta_r}$ on $V$.
\end{prop}
\begin{proof}
As we have seen in Section \ref{par:lattice}, the sections $e_i$ of $\cL_B(V)$ are locally of the form $[\diff
p_i]$, where each $p_i$ is a local Casimir on $B$. Thus, $X_{e_i}=X_{\phi^*p_i}$, so that the vector fields
$X_{e_i}$ are locally Hamiltonian vector fields, hence (globally) Poisson vector fields (on $\phi^{-1}(V)$). It
implies that for every function $H$ on $\phi^{-1}(V)$
$$ 
  [X_{e_i},X_{H}]=X_{X_{e_i}(H)}\;.
$$ 
In view of (\ref{eqn:angle_def}), this shows that $[X_{e_i},X_{\theta_j}]=0$ for $i,j=1,\dots,r$.  In turn, this
implies that for $F$ a function on $V$, the function $ X_{\theta_j}(\phi^* F) $ is a $ \phi$-basic function on
$\phi^{-1} (V)$; indeed, for any $i=1,\dots,r$,
$$ 
  X_{e_i}\left(X_{\theta_j}(\phi^* F)\right) =  X_{\theta_j}\left(X_{e_i}(\phi^*F)\right) =0\;.
$$
As a consequence, there exists a (unique) function $G_j$ such that $ \phi^* G_j = X_{\theta_j}(\phi^* F) $.  The
map $F \mapsto G_j $ is clearly a derivation, hence defines a vector field on $V$ which we denote by
$Y_{\theta_j}$.  By construction, the vector fields $X_{\theta_j} $ and $ Y_{\theta_j}$ are $ \phi$-related,
$\phi^*\circ Y_{\theta_j}=X_{\theta_j}\circ\phi^*$.  The fact that each $Y_{\theta_i}$ is a Poisson vector field
follows from the fact that $\phi $ is a Poisson submersion from $M$ to~$B$, and that $X_{\theta_i}$, which is a
Hamiltonian, hence Poisson vector field, is $\phi$-related to $Y_{\theta_i}$.  They commute in view of the
commutativity of the vector fields $X_{\theta_i}$ to which they are $\phi$-related, with $\phi$ being a submersion.
\end{proof}

We now show how two different sets of local angle variables are related.  Suppose that
$\theta=(\theta_1,\dots,\theta_r)$ and $\theta'=(\theta_1',\dots,\theta_r')$ are two sets of local angle variables
adapted to the same trivialization $(e_1,\dots,e_r)$ of $\cL_B(V)$.  Let $Y_{\theta_1}, \dots,Y_{\theta_r} $ be the
vector fields on $V$ defined in Proposition \ref{prop:angles_vectorfield_quotient} using the set of angle variables
$\theta_1,\dots,\theta_r$. Then there exist functions $F_1, \dots, F_r $ on $V$ such that:
\begin{enumerate}
  \item[(1)] $\theta_i' = \theta_i + \phi^* F_i $;
  \item[(2)] $ \pb{F_i,F_j} =  Y_{\theta_j}(F_i) -Y_{\theta_i}(F_j)$.
\end{enumerate}
Let us prove this claim. In view of (\ref{eqn:angle_def}), $ X_{e_j}(\theta_i - \theta_i')=0$ for all
$i,j=1,\dots,r$, which yields the existence of (unique) functions $F_1,\dots,F_r$ on $V$, satisfying (1). Since
$ \pb{\theta_i',\theta_j'}= \pb{\theta_i,\theta_j}=0$ for $i,j=1, \dots,r$, (1) implies:
$$ 
  0=\pb{\theta_i',\theta_j'}- \pb{\theta_i,\theta_j}= \pb{\theta_i,\phi^* F_j}+\pb{\phi^* F_i,\theta_j}+ 
  \pb{\phi^* F_i,\phi^* F_j}.
$$
Now, by definition of the vector fields $Y_{\theta_i}$ and since $\phi$ is a Poisson map, this amounts to:
$$ 
  \phi^* (Y_{\theta_i}(F_j))  -  \phi^*(Y_{\theta_j}(F_i)) -  \phi^* \pb{ F_i,F_j}  =0\;.
$$
This gives the second relation. Conversely, given a set of angle variables $\theta_1,\dots,\theta_r$ and functions
$F_1,\dots,F_r$ on $V$, satisfying (2), the above computation shows that the functions $\theta_i'$, defined by (1),
are a set of angle variables adapted to the same trivialization $(e_1,\dots,e_r)$ of $\cL_B(V)$.

\begin{remark}
For a given trivialization $e=(e_1,\dots,e_r)$ of $\cL_B(V)$, each one of the action variables $p_i$ satisfying
$e_i=[\diff p_i]$ is uniquely determined up to an element of $\Cas_B^M(V)$. Therefore, if a set of action variables
adapted to~$e$ exists, the space of all sets of action variables adapted to $e$ is an affine space of rank $r$ over
the ring $\Cas_B^M(V)$. There is no similar property for angle variables adapted to $ (e_1, \dots,e_r)$: it is not
an affine space, since the transformation which relates two of them (formulas (1) and (2) above) is non-linear.
\end{remark}

\subsection{Angle foliations}

For a given set of local angle variables $\theta=(\theta_1,\dots,\theta_r)$ on $\phi^{-1}(V)$, the level sets of
the map $\theta:\phi^{-1}(V)\to(\bbR/\bbZ)^r$ define a foliation $\fG_\theta$ of $\phi^{-1}(V)$, transverse to the
fibers of $\phi$, and having the following two properties:
\begin{enumerate}
  \item[(1)] $\fG_\theta$ is invariant under the flow of  the action vector field associated to any element of~$\cL_B(V)$;
  \item[(2)] $\fG_\theta$ is coisotropic, i.e., every leaf of $\fG_\theta$ is a coisotropic submanifold of $(M,\Pi)$.
\end{enumerate}
For the proof of (1), one needs to check that the Lie derivative with respect to the action vector fields $X_{e_i}$
of every first integral of $\fG_\theta$ is a first integral of~$\fG_\theta$; this is clear because the leaves of
$\fG_\theta$ are defined by $\theta_j=$ constant and $\cL_{X_{e_i}}(\theta_j)=X_{e_i}(\theta_j)$ is constant for
all $i$ and $j$, in view of (\ref{eqn:angle_def}). The proof of (2) follows from the fact that the functions
$\theta_j$, which define $\fG_\theta$, are in involution, again according to (\ref{eqn:angle_def}).

Making abstraction of these properties leads to the following definition.

\begin{defn}\label{def:angle_foliation}
Let $V$ be an open subset of $B$. A foliation $\fG$ of $\phi^{-1}(V)$ is called an {\bf angle foliation} if it has
the following properties:
\begin{enumerate}
  \item[(1)] $\fG$ is transverse to the fibers of $\phi$;
  \item[(2)] $\fG$ is invariant under the flow of  the action vector field associated to any element of~$\cL_B(V)$;
  \item[(3)] $\fG$ is coisotropic.
\end{enumerate}
\end{defn}

According to Proposition \ref{prp:angels_exist}, angle variables exist semi-locally, i.e., on an open neighborhood of any fiber of $\phi$, 
hence action foliations exist semi-locally. We show in the following proposition that every angle 
foliation is defined semi-locally by angle variables.

\begin{prop}\label{prop:from_angle_foliations_to_functions}
  Let $V$ be an open subset of $B$. We suppose that we are given on $V$ a trivialization $(e_1, \dots, e_r)$ of
  ${\cL_B(V)}$ and on $\phi^{-1}(V)$ an angle foliation $\fG$. Let $b\in V$. There exists a neighborhood $V_0$ of
  $b$, contained in~$V$, and there exist local angle variables $\theta=(\theta_1,\dots,\theta_r)$ on
  $\phi^{-1}(V_0)$, adapted to $(e_1, \dots, e_r)$, such that $\fG_\theta=\fG$ on $\phi^{-1}(V_0)$.
\end{prop}
\begin{proof}
It follows from (2) in Definition \ref{def:angle_foliation} that the flow of the (commuting) action vector fields
$X_{e_i}$ defines a diffeomorphism between $\phi^{-1}(V_0)$ and $\bbT^r\times V_0$ where $V_0$ is an open subset of
$V$ which contains $b$. By construction, this diffeomorphism has the following two properties: first, the
fundamental vector fields of the natural action of $\bbT^r$ on $\bbT^r\times V_0$ coincide with the action vector
fields $X_{e_i}$. Second, the leaves of $\fG$ correspond to the fibers of the projection map
$\theta:\phi^{-1}(V_0)\simeq\bbT^r\times V_0\to\bbT^r$; in particular, the foliations $\fG_\theta$ and $\fG$
coincide over points of $V_0$. Writing $\theta=(\theta_1,\dots,\theta_r)$ yields local angle coordinates on
$\phi^{-1}(V_0)$ adapted to $(e_1,\dots, e_r)$. Indeed, by construction, $X_{e_i} (\theta_j) =\delta_{i,j}$ for
$i,j=1, \dots,r$ and the functions $\theta_i$ are in involution because $\fG$ is coisotropic.
\end{proof}

The set of angle variables defining a given angle foliation is unique up to adding locally constant functions 
and taking integer-valued linear transformations. This is shown in  the following proposition. 

\begin{prop}\label{prop:unicity_up_to_constant}
Let $ \fG$ be an angle foliation on $\phi^{-1}(V)$, where $V$ is an open subset of $B$. Let $e=(e_1,\dots,e_r)$ and
$e'=(e_1',\dots,e_r')$ of $\cL_B(V) $ be two local trivializations of $V$ and denote by $C$ the invertible
integer-valued matrix such that $e'=e C$.  Let $\theta$ and $ \theta'$ be two sets of angle variables defining $\fG
$ and adapted to $e$ and $e'$ respectively. There exists a vector of locally constant $\bbR/\bbZ $-valued functions
$c=(c_1,\dots,c_r)$ on $\phi^{-1}(V)$, such that
\begin{equation}\label{eq:cocycleangle} 
  \theta'= \theta (C^t)^{-1}+c \;.
\end{equation}
\end{prop}
\begin{proof}
Suppose first that $e=e'$. Since both $\theta$ and $ \theta'$ define the same foliation~$\fG$, we have, in a
neighborhood of any point of $\phi^{-1}(V)$, $ \theta_i' = K_i(\theta_1, \dots,\theta_r)$ for some function $K_i$.
Applying $ X_{e_j}$ to both sides of the previous equation amounts to:
$$ 
  \delta_{i,j} = \sum_{k=1}^r \frac{\partial K_i}{\partial x_k} X_{e_j}(\theta_k)= \frac{\partial K_i}{\partial x_j}\; .
$$
This implies that $ \theta_i' - \theta_i$ is a locally constant function, which proves (\ref{eq:cocycleangle}) in
case $C=I_r$.  In general (i.e., without assuming that $e=e'$) the angle variables $ \theta' $ and $ \theta
(C^t)^{-1}$ are both adapted to $e'$, so that they differ by locally contant functions.
\end{proof}

The next theorem gives a necessary and sufficient condition for the existence of angle variables.  We use angle
foliations in its proof, in order to clarify the argument.

\begin{thm}\label{thm:angle_existence}
Let  $(M,\Pi) \stackrel{\phi}{\to }(B,\pi)$ be an NCI with compact fibers. The following statements are equivalent:
\begin{enumerate}
  \item[(i)] There exist global angle variables;
  \item[(ii)] The action lattice sheaf $\cL_B$ admits a global trivialization and there exists a section of $\phi :
    M \to B$ whose image is a coisotropic submanifold of $(M,\Pi)$.
\end{enumerate} 
\end{thm}
\begin{proof}
As pointed out after Definition \ref{def:angles_variables}, if there exists a set of global angle variables
$(\theta_1, \dots,\theta_r)$, then the action lattice sheaf $\cL_B$ admits a global trivialization.  The zero locus
$ \theta_1 = \dots =\theta_r=0$ is a submanifold $B_0$ which is transverse to the fibers of $\phi : M \to B$. Since
the restriction of $\phi$ is a bijection from $B_0$ to $B$, it is the image of some section $\sigma$ of $ \phi: M
\to B$. Since the foliation $\fG_\theta$ which is associated to $\theta$ is coisotropic, $B_0$ is coisotropic. This
proves $(i)\implies (ii)$.

Let us prove that $(ii)$ implies $(i)$. A choice of global trivialization $(e_1, \dots,e_r)$ of $\cL_B$ turns $M
\to B$ into a principal $\bbT^r$-bundle; we denote by $(s,m) \to s \cdot m$ the action of $s \in \bbT^r$ on $m \in
M$. Let $\sigma : B \to M$ be a section of $\phi : M \to B$ whose image $B_0:=\sigma(B)$ is a coisotropic
submanifold. Consider the unique $\bbT^r$-invariant foliation $\fG$ on $M$ admitting $B_0$ as a leaf, i.e. consider
the foliation admitting the submanifolds $ s\cdot B_0$ with $s \in \bbT^r$ as leaves. By construction, $\fG$ is
transverse to all fibers of $\phi$. Also, $\fG$ is $\bbT^r$-invariant, so that it is invariant under all the action
fields associated to elements of~$\cL_B(B)$. Since for all $ s \in \bbT^r$, the map $m \to s\cdot m $ is a Poisson
diffeomorphism of $M$, the fact that $B_0$ is a coisotropic submanifold implies that all the leaves of $\fG$ are
coisotropic submanifolds, so that $\fG$ is an angle foliation.

According to Proposition \ref{prop:from_angle_foliations_to_functions}, there exists for any $b \in B$ a
neighborhood $U_b$ of $\phi^{-1}(b)$ and a unique set $(\theta_1, \dots,\theta_r)$ of angle variables on $U_b$,
adapted to $(e_1, \dots,e_r)$, constant on the leaves of $\fG $ and vanishing on $B_0$.  The open subsets $
(U_b)_{b \in B}$ form an open cover of $M$. Since the angle variables defined on $U_b$ and $U_b'$ coincide on $U_b
\cap U_b'$, they lead to global angle variables.
\end{proof}

\subsection{The transverse foliation}\label{par:transverse_foliation}
We have seen in Section \ref{par:angle_variables} that we can associate to a set of local angle variables
$\theta=(\theta_1,\dots,\theta_r)$ on $\phi^{-1}(V)$ vector fields $Y_{\theta_1},\dots,Y_{\theta_r}$ on $V\subset
B$. We now show that they define a distribution of rank $r$ on $V$ which depends only on the angle foliation,
defined by the angle variables. For a given set of local angle variables, let us denote by $D_\theta$ the (a priori
singular) distribution on $V$, defined by the vector fields $Y_{\theta_1},\dots,Y_{\theta_r}$, where
$Y_{\theta_i}:=\phi_*X_{\theta_i}$ and by $L_\theta$ the (a priori singular) lattice subbundle of~$D_\theta$,
generated by these vector fields.
\begin{prop}\label{prop:angles_distribution_quotient}
  Let $V$ be an open subset of $B$ and suppose that $\fG$ is an angle foliation on $\phi^{-1}(V)$, where $V$ is an
  open subset of $B$. Suppose that $\fG$ is defined by local angle variables $\theta=(\theta_1,\dots,\theta_r)$.  
\begin{enumerate}
  \item[(1)] $D_\theta$ is an integrable distribution of rank $r$ on $V$;
  \item[(2)] $D_\theta$ and $L_\theta$ are independent of the choice of $\theta$, defining $\fG$.
\end{enumerate}
Therefore, $\fG$ defines an $r$-dimensional foliation $\fT_\fG$ of $V$ and a lattice bundle~$L_\fG$ on~$V$, which
we call the {\bf transverse foliation}, respectively the {\bf transverse lattice bundle} of the NCI system.
\end{prop}
\begin{proof}
Using the angle foliation $\fG$ we can define an $r$-dimensional subspace $D'_m$ of $T_mM$ at very point $m\in
\phi^{-1}(V)$ by setting $D_m':=\Pi_m^\sharp((T_m\fG)^0)$. It leads to a distribution $D'$ on $\phi^{-1}(V)$, which
is spanned by the $r$ independent commuting vector fields $X_{\theta_i}$ at $m$, where
$\theta=(\theta_1,\dots,\theta_r)$ is any set of local angle variables defining $\fG$ around $m$. Thus, its
projection under $\phi$, whose fibers are transverse to $\fG$, is a distribution which is spanned by the $r$ vector
fields $Y_{\theta_i}$ on $B$, hence it is the distribution $D_\theta$.  It follows that $D_\theta$ is an integrable
distribution of rank $r$ on $V$ and that $D_\theta$ is independent of the choice of $\theta$, defining $\fG$. The
integral manifolds of $D_\theta$ are the leaves of an $r$-dimensional foliation of $V$, denoted by $\fT_\fG$. In
view of (\ref{eq:cocycleangle}), two different choices $\theta$ and $\theta'$ are related by $\theta'= \theta
(C^t)^{-1}+c,$ where $C$ is an integer-valued matrix and $c$ is a constant vector. It follows that $L_\theta$ and
$L_{\theta'}$ define the same lattice bundle in $D_\theta=D_{\theta'}$.
\end{proof}

We show in the following proposition how an action and an angle foliation, if they exist, are related.
\begin{prop}\label{prop:angle_foliation_distribution_quotient}
  Let $V$ be an open subset of $B$. Suppose that we have on~$V$ an action foliation $\fA$ and on $\phi^{-1}(V) $
  an angle foliation $\fG$. 
  \begin{enumerate} 
    \item[(1)] $\fT_\fG $ is transverse to $\fA $;
    \item[(2)] The tangent space to $\fT_\fG $ is generated by (local) Poisson vector fields which preserve $\fA$.  
  \end{enumerate}
\end{prop}
\begin{proof}
In a neighborhood $V_0$ of any point of $V$, there exist action-angle variables
$p_1,\dots,p_r,\theta_1,\dots,\theta_r$ such that $([\diff p_1],\dots,[\diff p_r])$ is a trivialization of
$\cL_B(V_0)$. Hence:
\begin{equation}\label{eq:Y_i_And_p_i} 
  \phi^*(Y_{\theta_i}(p_j)) = X_{\theta_i} (\phi^*p_j) = - X_{\phi^* p_j } (\theta_i) = -X_{e_j}( \theta_i )= -\delta_{i,j}\;,
\end{equation} 
which implies both items (1) and (2).
\end{proof}

Consider a foliation $\fG $ of $\phi^{-1}(V)$ transverse to the fibers of the surjective submersion $\phi : M \to
B$, where $V$ is an open subset of $B$. For any leaf $G$ of~$\fG $, $ \phi$ is a local diffeomorphism from $G$ to
$B$, so that a multivector field on $B$ induces a multivector field on the leaf $G$. Making this construction for
all the leaves of $\fG$ simultaneously, yields a graded Lie algebra morphism $ \phi_\fG^*$ from the space of
multivector fields on~$V$ to the space of multivector fields on $\phi^{-1}(V)$ tangent to the foliation $\fG$,
where both spaces are equipped with the Schouten bracket.

We apply this to the case of an angle foliation $\fG$ on $\phi^{-1}(V)$ with $V$ an open subset of $B$, to
construct two Poisson structures on $\phi^{-1}(V)$, to wit~$\phi_\fG^* (\pi) $ (with $\pi$ the Poisson structure on
$B$) and
\begin{equation}\label{eq:bivecteur_regulier} 
  \Pi_\fG := \sum_{i=1}^r X_{e_i} \wedge \phi_\fG^* (Y_{\theta_i})\; .
\end{equation}
In this formula, the $\theta_i$ stand for any set of local action variables, defined in a neighborhood $W$ of some
point of $V$ and $e=(e_1, \dots,e_r)$ stands for the corresponding trivialization of $\L_B(W)$ and $
Y_{\theta_1},\dots,Y_{\theta_r}$ are the vector fields on $W$, defined in Proposition
\ref{prop:angles_vectorfield_quotient}; the right hand side of (\ref{eq:bivecteur_regulier}) does not depend on the
choice of $\theta_i$ because the $\theta_i$, and hence the vector fields $Y_{\theta_i}$, are dual to the
trivialization $e$. It follows that the right hand side of (\ref{eq:bivecteur_regulier}) is a well-defined bivector
field on~$\phi^{-1}(V)$.

\begin{prop}\label{prp:decomp}
Let $V$ be an open subset of $B$ and suppose that $\fG$ is an angle foliation on $\phi^{-1}(V)$. 
\begin{enumerate}
  \item[(1)] The bivector field $\Pi_\fG $ is a regular Poisson structure on $\phi^{-1}(V)$ of rank $2r$.
  \item[(2)] The Poisson structures $ \Pi$, $ \Pi_\fG$ and $\phi_\fG^* (\pi) $ are related by:
          $$ \Pi = \Pi_{\fG} + \phi_\fG^*(\pi). $$
\end{enumerate}
\end{prop}
\begin{proof}
Let us first rewrite the local expression of $\Pi_\fG $ given in formula (\ref{eq:bivecteur_regulier}) in a more
convenient way. Choose a trivialization $e=(e_1,\dots,e_r) $ of $\cL_B(V) $, a set of local angle variables
$(\theta_1,\dots,\theta_r)$ adapted to $e$ defining $\fG$, and a set of local action variables $p=(p_1, \dots,p_r)
$ satisfying $ e_i =[\diff p_i]$.  For $i=1, \dots,r$, the identity $X_{\phi^* p_i}=X_{e_i} $ holds. Also,
$Y_{\theta_i}$ is $\phi$-related to $ X_{\theta_i}$, which is tangent to $\fG$, so that $ X_{\theta_i}
=\phi_\fG^*(Y_{\theta_i})$. It follows that (\ref{eq:bivecteur_regulier}) can be written as
\begin{equation}\label{eq:explicit_piG} 
  \Pi_\fG = \sum_{i=1}^r  X_{\phi^* p_i} \wedge X_{\theta_i}\;.
\end{equation}
Since the $2r$ vector fields $ X_{\phi^* p_1},\dots, X_{\phi^* p_r} , X_{\theta_1}, \dots , X_{\theta_r}$ are
pairwise commuting, $\Pi_\fG$ is a Poisson structure.  Also, (\ref{eq:explicit_piG}) implies that $\Pi_\fG (\diff
p_j,\diff \theta_i) = \delta_{i,j}$ while $\Pi_\fG (\diff \theta_i,\diff\theta_j) = \Pi_\fG (\diff p_i,\diff p_j) =
0 $, which proves that $\Pi_\fG $ is a regular bivector field of rank $2r$. This proves~(1).

The bivector field $ P:=\Pi - \Pi_\fG $ is tangent to $\fG $, i.e. $ P_m \in \wedge^2 T_m \fG$ for every $m
\in\phi^{-1}(V)$. Indeed, we have in view of (\ref{eq:explicit_piG}) that $\Pi_\fG^\sharp(\diff
\theta_j)=-X_{\theta_j}=\Pi^\sharp(\diff\theta_j)$.  Also, $\wedge^2 T_m \phi (\Pi_\fG)_m=0$ so that $ \wedge^2
T_m\phi (\Pi_m) = \pi_{\phi(m)}$. This shows that on $\phi^{-1}(V)$ both bivector fields $P$ and $ \phi_\fG^* (\pi)
$ are tangent to $\fG$ and project to $\pi $, so they are equal and (2) follows.
\end{proof}

The difference between the existence of angle foliations and angle variables can also be stated in the following
geometrical terms. Suppose that $(M,\Pi) \stackrel{\phi}{\to }(B,\pi) $ is an NCI system with {compact fibers} and
suppose that its lattice sheaf $\cL_B$ admits a global trivialisation, so that $M\to B$ is a principal
$\bbT^r$-bundle. The distribution, tangent to an angle foliation $\fG$ is an Ehresmann connection, which is
invariant under the torus action, hence it defines a principal $\bbT^r$-connection. By construction, this
distribution is integrable, which is tantamount to saying that the connection is flat. Saying that there exist
angle variables, adapted to $\fG$ is equivalent to saying that the bundle $M\to B$ is trivial, hence is of the
form $\phi:\bbT^r\times B\to B$, where $\phi$ is the projection on the second component and the $\bbT^r$ action is
the standard one.

\subsection{The transverse Poisson manifold}
In this paragraph we give necessary and sufficient conditions for $(M,\Pi)$ to be Poisson diffeomorphic with the
product $\bbT^r\times T\times A$, where $A$ is a leaf of $\fA$, equipped with the Poisson structure inherited from
$(B,\pi)$ (as a Poisson submanifold) and $T$ is a leaf of~$ \fT_\fG$, the Poisson structure on $\bbT^r\times W$
being the canonical Poisson structure defined by a set of global action-angle variables, which we assume to exist. 

In order to do this, we first recall a basic result from foliation theory. Suppose that $\fA$ and $\fT$ are two
foliations of a manifold $B$ which intersect transversally (as the notations suggest, we will use the result when
$\fA$ and~$\fT$ are the action and transverse foliations on $B$, defined by the action-angle variables). We say
that $\fA$ and $\fT$ have the {\bf unique intersection property} if any leaf of $\fA$ has exactly one point in
common with any leaf of $\fT$. Fix a point $b\in B$ and denote by $A$ and $T$ the leaves of $\fA$ resp.\ of $\fT$,
passing through~$b$.  There is a neighborhood $ V_b$ of $b$ in $B$ and a unique diffeomorphism $\Phi_b$
from~$V_b$ to $A_b \times T_b$ with $A_b$ and $T_b$ a neighborhood of $b$ in $A$ resp.\ in $T$, under which the
foliations $\fA$ and $\fT$ become the fibers of the projections onto the first and second components
respectively. Since this diffeomorphism on~$V_b$ is unique, it leads to a global diffeomorphism between $B$ and
$A\times T$ if (and only if) the foliations $\fA$ and $\fT$ of $B$  have the unique intersection property.

\begin{thm}\label{thm:aa_when_exists}
  Let $(M,\Pi) \stackrel{\phi}{\to}(B,\pi)$ be a NCI system with compact fibers, equipped with a set of angle
  variables $\theta := (\theta_1,\dots,\theta_r) $ and a set of action variables $p := (p_1, \dots,p_r) $.  We set
  $W:=p(B) $, which is a connected open subset of ${\mathbb R}^r $. Choose a point $b \in B$ and let $A$ and $T$
  denote the leaves through $ b$ of the action foliation $\fA$, associated to $p$ and of the transverse foliation
  $\fT_\fG$, associated to $\theta$. Then the following are equivalent:
  \begin{enumerate}
    \item[(i)] The map $p$ restricts to a bijection from $T$ to $W$, and the foliations $\fA $ and $ \fT_\fG$ have
      the unique intersection property.
    \item[(ii)] There exist diffeomorphisms $\chi$ and $\chi_B$ making the following diagram commutative:
    $$
      \xymatrix{
        &M\ar[r]^{\phi}\ar[dl]_{\theta}\ar[dd]^{\chi}_{\simeq}&B\ar[dd]^{\chi_B}_{\simeq}\ar[rd]^{p}\\  \bbT^r&&&W\\
         &\bbT^r\times T\times A\ar[ul]\ar[r]&T\times A\ar[ur]}
    $$ 
    \end{enumerate}
  Moreover, when these conditions are satisfied, 
  \begin{equation}\label{eq:poisson_push}
    \chi_* (\Pi) = \sum_{i=1}^r \frac{\partial}{\partial \theta_i} \wedge \frac{\partial}{\partial p_i} + \pi_{\vert_A}\ 
    \quad\hbox{and}\quad (\chi_B)_*\pi=\pi_{\vert_A}\;. 
  \end{equation}
\end{thm}
\begin{proof}
Recall that the action and transverse foliations, when they exist, are transverse. We assume here to be given
global action-angle variables, hence both foliations exist and we can apply the above remarks on transversally
intersecting foliations to prove the equivalence of (i) with the existence of $\chi_B$ in (ii), making the
rightmost triangle in the above diagram commutative. In view of the existence of action-angle variables, $M$ is a
trivial $\bbT^r$-bundle over $B$, allowing us to complete the diagram. This shows the equivalence of (i) and (ii).

Locally, $\chi_B$ is a a Poisson diffeomorphism between an open neighborhood in $B$ and open neighborhoods in the
leaves $A$ and $T$, when $A\times T$ is equipped with the product of $\pi$ restricted to the Poisson submanifold
$A$ and the trivial Poisson structure on $T$. This follows from the fact that the foliation $\fT_\fG$ is generated
by Poisson vector fields which preserve the foliation $ \fA$ (see Propositions
\ref{prop:angles_distribution_quotient} and \ref{prop:angle_foliation_distribution_quotient}). Since $\chi_B$ is a
(global) diffeomorphism, it is a Poisson diffeomorphism, leading to the second formula in
(\ref{eq:poisson_push}). The first formula in (\ref{eq:poisson_push}) follows from Proposition \ref{prp:decomp},
\end{proof}


\section{Examples}\label{section:examples}

In this section we give a series of examples and counter-examples which illustrate the different obstructions to
the existence of global action-angle variables.

\subsection{An isotropic Poisson complete foliation which is not an abstract NCI system}\label{ex:not_NCI}

We first give an example which shows that not every Poisson complete foliation is an abstract NCI system. Consider
the trivial circle bundle $M:= S^1\times\bbR^3\to\bbR^3$ over~$\bbR^3$. Denoting the coordinates on $S^1$ and on
$\bbR^3$ by $\theta$ and $x,y,z$ respectively, we consider on $M$ the Poisson structure
\begin{equation*}
  \Pi:=\frac{\partial}{\partial \theta} \wedge \frac{\partial}{\partial z}+ \pi\;,
\end{equation*}%
where $\pi$ is the Poisson structure on $\bbR^3$ (or on $M$), given by
\begin{equation*}
  \pi:=\left(y\pp{}{x} - x\pp{}{y}\right)\we\pp{}z+(x^2 +y^2)\pp{}{x}\wedge\pp{}{y}\;.
\end{equation*}
Using the fact that $\left(y\pp{}{x} - x\pp{}{y}\right)(x^2 +y^2)=0$, one easily checks that $\pi$ and $\Pi$ are
indeed Poisson structures. Also, by construction, the canonical projection $\phi:(M,\Pi)\to(\bbR^3,\pi)$ is a
Poisson map. According to Example \ref{exa:poisson_submersion}, the fibers of $\phi$, which are circles, define a
Poisson complete foliation $\fF$ of $(M,\Pi)$. To see that $\fF$ is isotropic, take a point
$m=(\theta_0,x_0,y_0,z_0)$ of $M$ and consider $\a_m=a\,\diff x+b\,\diff y-\diff z$, where $a,b\in\bbR$. By a
direct computation we find that $\Pi^\sharp(\a_m)=\p/\p\theta+\pi^\sharp(\a_m)=\p/\p\theta$ when $a$ and $b$ are
taken as
\begin{equation*}
  a=\frac{x_0}{x_0^2+y_0^2}\;,\qquad b=\frac{y_0}{x_0^2+y_0^2}\;;
\end{equation*}%
for $x_0=y_0=0$ these formulas do not make sense, but in that case any values of $a$ and $b$ do the job. Since
clearly $\a_m\in (T_m\fF)^\circ$, this shows that~$\fF$ is isotropic. We now show that in a neighborhood $U$ of
$m=(\theta_0,0,0,z_0)$ there exists no function $f$, constant on the leaves of $\fF$, whose Hamiltonian vector
field $\X_f$ generates $T\fF$ on $U$. The first condition means that~$f$ is independent of~$\theta$, so that
\begin{equation}\label{eq:ders_co_ex}
  \diff \theta(X_f)=\pp fz\;,\qquad \diff x(X_f)=y\,\pp fz+(x^2+y^2)\pp fy\;.
\end{equation}%
The second condition means that $X_f=g\,\p/\p\theta$, for some nowhere vanishing function $g$ on $U$, so that
$\diff \theta(X_f)\neq0$ and $\diff x(X_f)=0$ on $U$. In view of (\ref{eq:ders_co_ex}) this is impossible.

\subsection{The existence of action variables and foliations}
We now give two examples of NCI systems which have compact fibers and trivial action lattice sheaf, yet fail to
have action variables; the two examples differ in the existence of a global action foliation. We also show that the
existence of action variables, defining an action foliation may depend on the choice of action foliation.

Let $ M:= S^1 \times B$ where $B$ is a manifold equipped with a nowhere vanishing vector field $\cV$. The foliation
of $B$, defined by $\cV$, is denoted by $\fT$. Consider the Poisson structure on $M$ defined by
$$ 
  \Pi := \pp{}{\theta} \we \cV\;,
$$
where $\theta $ is the parameter on $S^1$, viewed as a function on $M$. Let $\phi:S^1\times B\to B$ denote the
projection on the second component. The tangent space to the fibers of $\phi$ is generated $\p/\p\theta$, which is
a locally Hamiltonian vector field: for any local function $p$ on~$B$ we have that
$\X_{\phi^*p}=\phi^*(\cV(p))\,\p/\p\theta$.  Thus, $(M,\Pi)\stackrel{\phi}{\to}(B,0)$ is an NCI system of rank
1. The fibers of its momentum map are circles. For every point $b \in B$, only one of the two generators of the
action lattice $L_b$ at the point $b$ corresponds to the vector field $\pp{}{\theta}$. The action lattice,
therefore, admits a global section $e$, in particular the action lattice sheaf $\cL_B$ is trivial.
\begin{prop} \label{prop:no_action_foliation}
  When $B$ is compact, the NCI system $(M,\Pi)\stackrel{\phi}{\to} (B,0)$ above does not admit global action
  variables.  When $B$ is moreover simply-connected, it even does not admit an action foliation.
\end{prop}
\begin{proof}
When $B$ is compact, every function on $B$ has points where its differential vanishes. Such a function can never be
an action variable, which shows the first statement. Assume now that there exists a global action foliation $\fA$
on $B$. By passing to the orientation cover, we can assume that $\fA$ is co-oriented. 
Since the rank of the NCI system is $1$, $\fA$ is a transverse integral affine foliation of codimension $1$, so it must 
be given by the kernel of a \emph{closed} 1-form. When $B$ is simply-connected, $H^1(B,\bbR)=0$, so this form is exact
and its kernel cannot define a regular foliation.
\end{proof}
The second part of this proof can be reformulated in terms of the obstruction theory of 
Section \ref{sec:action:foliation} as follows: according to Proposition
\ref{cor:actionfoliation}, an action foliation exists iff $\Obs([e])$ is representable by constants. When $B$ is simply-connected,
$H^1(B,\bbR)=0$, so $\Obs([e])$ is representable by constants if and only if $\Obs([e])= 0$, which is according to
Theorem \ref{prp:action_vars_exist} equivalent to the existence of a global action variable. But we know, from the
first part that such a global variable does not exist. 

Let us apply the proposition to $B=S^3$, equipped with the fundamental vector field $\cV$ of the Hopf fibration
$S^3 \to S^2 $, {i.e.,} the fundamental vector field of the natural $S^1$-action on $S^3$. Since $S^3$ is both
compact and simply-connected, Proposition \ref{prop:no_action_foliation} shows that this NCI system that does not
admit an action foliation.

We next apply the proposition to $B=S^1$, with its natural translation invariant vector field $\p/\p\psi$, so that
$\omega=\diff\psi$. Since $S^1$ is compact, Proposition \ref{prop:no_action_foliation} shows that this system does
not admit an action variable. However, since $\omega$ is closed (but not exact!), it defines an action foliation.

To finish, we consider $B:=S^1\times \bbR$ (a cylinder) equipped with an $S^1$-valued coordinate $\psi$ and an
${\bbR}$-valued coordinate $p$, corresponding to the first and second projections. Any foliation $\fA$ of $B$,
transverse to $\cV:=\p/\p p$ is an action foliation since $\fA$ can locally be defined by a function $\tilde{p}$
such that $\pp{\tilde{p}}p=1$, i.e., a local action variable. Thus the two foliations,  defined by the vector
fields 
$$ 
\pp{}\psi  \qquad\hbox{and}\qquad \pp{}\psi + p \pp{} p \;
$$
are action foliations. The first foliation is defined by the function $p$, which is an action variable.
However, the second foliation has as leaves the circle $C_0:=\{p=0\} $ and a family of curves which are transverse
to $\p/\p p$ and spiral towards $C_0$.  It is not a foliation defined by a function, so there is no global action
variable defining it.

\subsection{The existence of  angle variables and foliations}
Consider an NCI system $(M,\Pi)\stackrel{\phi}{\to}(B,\pi)$ of rank $r=1$ with compact fibers. We assume that its
action lattice sheaf admits a trivialization. Recall from Remark \ref{rem:LB_trivial} that this implies that
$\phi:M{\to}B$ is a principal $S^1$-bundle. Notice that in the rank 1 case every section of $\phi:M\to B$ is
coisotropic, because the image of such a section is of codimension 1. It follows that the principal $S^1$-bundle
$\phi:M{\to}B$ has the following properties:
\begin{enumerate}
  \item[(1)] It admits a trivialization if and only if there exists a global angle variable;
  \item[(2)] It admits a flat connection if and only if there exists a global angle foliation.
\end{enumerate}

Indeed, Theorem \ref{thm:angle_existence} yields in the present case that a global angle variable exists if and
only if a global section of $\phi$ exists, which is itself equivalent to the triviality of the principal
$S^1$-bundle. This shows (1). Also, the connection form of a principal $S^1$-bundle is simply a nowhere vanishing
one-form $\beta \in \Omega^1(M,{\bbR})$, and such a connection is flat if and only if $\beta$ is closed, which in
turn implies that the distribution $\Ker\b$ is integrable, hence defines a foliation transverse to the fibers of
$\phi$. It is an angle foliation, because it is of codimension 1 (hence coisotropic) and because the connection
form $\beta$ is $S^1$-invariant.  Conversely, the leaves of any angle foliation of the NCI system define an
integrable distribution which is transverse to the fibers of $\phi$ and is $S^1$-invariant, i.e. a flat
connection. This shows (2).

Let $\phi_0:M_0{\to}B_0$ be a principal $S^1$-bundle and denote the fundamental vector field of the $S^1$-action on
$M_0$ by $\cW$. We associate to it an NCI system $ (M,\Pi) \stackrel{\phi}{\to} (B ,0)$ of rank 1 by setting
$M:=M_0\times \bbR$, $B:=B_0\times\bbR$ and $\phi:=\phi_0\times\hbox{Id}_\bbR$. The Poisson structure on $M$ is given
by $\Pi:=\pp{}{p}\we\cW $, where $ p$ is the parameter on ${\bbR}$. Clearly, the NCI system has compact fibers and
its action lattice sheaf admits a trivialization; indeed, $\phi:M{\to} B$ is a principal $S^1$-bundle. This bundle
admits a flat connection (respectively, is trivial) if and only if $\phi_0:M_0{\to}B_0$ admits a flat connection
(respectively, is trivial). Therefore, in order to construct an NCI system with compact fibers which admits no
angle foliation and an NCI system with compact fibers that admits an angle foliation but no angle variables, it
suffices to find:
\begin{enumerate}
  \item[(A)] A principal $S^1$-bundle which does not admit a flat connection;
  \item[(B)] A non-trivial principal $S^1$-bundle which admits a flat connection.
\end{enumerate}
The Hopf fibration $S^3\to S^2$ is an example of (A). In order to give an example of (B) we consider on $S^2\times
S^1$ the equivalence relation $R$ defined by $(x,y)\sim (-x,-y)$. The quotient map $S^2\to\bbR\bbP^2$ leads to a
map $\phi_0:(S^2\times S^1)/R\to\bbR\bbP^2$ which makes it into a non-trivial principal $S^1$-bundle. The standard
vector field $\p/\p\theta$ on $S^1$ is invariant under $y\mapsto -y$, hence leads to a non non-vanishing vector
field on $(S^2\times S^1)/R$ which is both $S^1$-invariant and transverse to the fibers of $\phi_0$. It defines a
distribution on $(S^2\times S^1)/R$ which is a flat connection.

As in the case of action variables, an NCI system may have two different angle foliations, where one can be defined by
angle variables while the other one can't. In view of the above analysis, an example for $r=1$ can be constructed
from a trivial  $S^1$-bundle $M=S^1\times B$ with two flat connections, one which is associated to a trivialization
but not the other one. We can take $B:=S^1$ and choose for the second connection and translation invariant
distribution on the torus $M$ whose leaves spiral at least twice around the torus.

\subsection{Sections versus coisotropic sections of the momentum map}
We have seen in Theorem \ref{thm:angle_existence} that global angle variables can only exist when the momentum map
has a coisotropic section. We now show that a coisotropic section of the momentum map may fail to exist even when
the momentum map has a section. Our example admits both an action foliation and a trivialization of its action
lattice sheaf.

We consider the NCI system $ (M,\Pi) \stackrel{\phi}{\to} (B ,0)$ where $M:=\bbT^2\times B$, where $\phi$ is the
projection on the second component and $B:=\bbT^2$. Also,  $\Pi$ is given by
 \begin{equation*}
 \Pi:=\pp{}{\theta_1}\we\pp{}{\psi_1}+\pp{}{\theta_2}\we\pp{}{\psi_2}+\alpha \pp{}{\theta_1}\we\pp{}{\theta_2}.
 \end{equation*}%
where $\alpha\in\bbR^*$, the standard ($S^1$-valued) coordinates on $B$ are denoted by $(\psi_1,\psi_2)$ and those
on the first factor of $M$ by $(\theta_1,\theta_2)$. Throughout the example we identify $S^1$ with $\bbR/\bbZ$ and
$\bbT^2$ with $S^1 \times S^1 $. The action lattice sheaf $\cL_B$ admits $(e_1,e_2):=([\diff\psi_1],[\diff\psi_2])$
as trivialization and we have $X_{e_i}(\theta_j)=\delta_{i,j}$. However, $(\theta_1,\theta_2)$ is not a set of
angle variables because $\pb{\theta_1,\theta_2}=\alpha$. If $(\theta'_1,\theta'_2)$ is a set of angle variables
adapted to the trivialization $(e_1,e_2)$, then $\theta_i'=\theta_i+\phi^*F_i$, for some $ S^1$-valued functions
$F_1,F_2$ on $B$; also, if we want that $\theta_1'=\theta_2'=0$ defines a coisotropic submanifold, we must have
$\pb{\theta_1',\theta_2'}=0$, to wit
\begin{equation}\label{eq:impossible}
  \alpha-\pp{F_1}{\psi_2}+\pp{F_2}{\psi_1}=0\;.
\end{equation}
Let $F$ be any smooth map from $S^1=\bbR/\bbZ$ to itself. Since any two smooth liftings $\tilde F:\bbR\to\bbR$
differ by an integer, the integral $\int_{S^1}F\diff\psi$ is well-defined up to an integer and $\int_{S^1}
\pp{ F}{\psi} \diff \psi \in {\bbZ} $. Therefore,
$$ 
  \int_{ S^1} \pp{F_1}{\psi_2} \diff \psi_2 \in {\bbZ} \quad\hbox{and}\quad
  \int_{ S^1} \pp{F_2}{\psi_1} \diff \psi_1 \in {\bbZ} \;,
$$
so that
$$ 
  \iint_{ S^1\times S^1}\big(\pp{F_1}{\psi_2}-\pp{F_2}{\psi_1} \big) \;\diff \psi_1 \diff\psi_2 \in {\bbZ}\;.
$$
However, $\iint_{S^1\times S^1} \alpha \diff \psi_1 \diff \psi_2 =\alpha$, so there is no solution to Equation
(\ref{eq:impossible}) unless $ \alpha\in\bbZ$. This shows that a set of angle variables adapted to the
trivialization $(e_1,e_2)$ does not exist, hence no set of angle variables exists (see Proposition
\ref{prop:unicity_up_to_constant}). In, turn, this implies that no coisotropic section of the momentum map of this
NCI system exists.

\subsection{The Euler-Poinsot top}

The configuration space of the Euler-Poinsot top is the Lie group $\G:={\rm SO}(3)$ of real orthogonal $3\times 3$
matrices, so its phase space is the cotangent bundle $T^*\G$, equipped with its canoncial symplectic
structure. Denoting the Lie algebra of $\G$ by $\fg$, we have that $T^*\G\simeq \G\times\fg^*$, where the
isomorphism is constructed by using left translation on $\G$. It is well-known that the symplectic manifold
$\G\times\fg^*$ is a symplectic groupoid in the sense of \cite{CDW}, with target map $t:\G\times\fg^*\to\fg^*$ the
(coadjoint) action map $(g,\xi) \mapsto \Ad_g^* \xi$ and source map $s:\G\times\fg^*\to\fg^*$ the projection onto
the second component, $(g,\xi) \mapsto\xi$. Like for any symplectic groupoid,
\begin{itemize}
  \item The source map $s$ is a Poisson map onto $\fg^*$, equipped with its Lie-Poisson structure;
  \item The target map $t$ is an anti-Poisson map onto the same space;
  \item For every pair of functions $F,G$ on $\fg^*$, the functions $s^*F$
and $t^*G $ are in involution on $\G\times \fg^*$.
\end{itemize}
It is convenient to identify $\fg^*$ with $\bbR^3$. First, we can identify $\fg^*$ with $\fg$ by using the Killing
form. Next, $\fg$ is the Lie algebra $\mathfrak{so}(3)$ of real skew-symmetric $3\times 3$ matrices, which we can
identify with $\bbR^3$ by assigning to $(x,y,z)\in\bbR^3$ the skew-symmetric matrix
$
  \begin{pmatrix}
   0 &z & -x\\  -z&0 & y \\ z& -y &0
  \end{pmatrix}.
$
Under these identifications:
\begin{itemize}
  \item The coadjoint action of ${\rm SO}(3)$ on $\mathfrak{so}(3)^*$ becomes the canonical action of ${\rm SO}(3)$
    on $\bbR^3$;
  \item The Lie bracket on $\fg$ becomes the vector product on $\bbR^3$;
  \item The Lie-Poisson structure on $\fg^*$ becomes the linear Poisson structure on $\bbR^3$, given in terms of
    the natural coordinates $(x,y,z)$ on ${\bbR}^3$ by:
    \begin{equation}\label{eq:LiePoisson}
      \pb{x,y}_{\fg^*}=z \,,\quad \pb{y,z}_{\fg^*}=x\,,\quad\pb{z,x}_{\fg^*}=y\;.
    \end{equation}
    A Casimir of this Poisson structure is given by $C:=x^2+y^2+z^2$.
\end{itemize}
The upshot is that ${\rm SO}(3) \times {\bbR}^3 $ is a symplectic manifold, comes equipped with two maps $s,t:{\rm
  SO}(3) \times {\bbR}^3\to\bbR^3$ which are defined by $s(R,m)=m$ and $t(R,m)=Rm$ and which are Poisson,
resp.\ anti-Poisson maps. Also, for every pair of functions $F,G$ on ${\bbR}$, the functions $s^*F$ and $t^*G $ are
in involution.  In turn, this implies that for any function $H$ on ${\bbR}^3 $, the map $\phi_H$, defined by
$$
  \begin{array}{lcccl}
    \phi_H &:& {\rm SO}(3)\times\bbR^3 &\mapsto & {\bbR}^3 \times {\bbR} \\
            && (R,m)                   &\to     & \big( Rm,H(m) \big).
    \end{array}
$$
is a Poisson map, when $ {\bbR}^3 \times {\bbR} $ is equipped with the Poisson structure $\pi=\PB_B$, which is the
product of the linear Poisson structure (\ref{eq:LiePoisson}) on $ {\bbR}^3 $ with the trivial Poisson structure on
${\bbR}$. The symplectic Poisson structure on ${\rm SO}(3) \times {\bbR}^3$ is denoted by $\Pi=\PB$.

The \emph{Euler-Poinsot top} corresponds to the choice
$$
  H:= \frac{1}{2} \left( \frac{x^2}{I_x} + \frac{y^2}{I_y}+\frac{z^2}{I_z}\right)\;.
$$
where $I_x\,I_y$ and $I_z$ are positive parameters, describing the top. In what follows we assume that these
parameters are different and that the coordinates are ordered such that $I_x>I_y>I_z$. Consider the functions
$s^*H,t^*C,t^*x,t^*y$ and $t^*z$ on ${\rm SO}(3) \times {\bbR}^3$ and consider the Hamiltonian vector fields
$\X_{s^*H}$ and~$\X_{t^*C}$. On the one hand, $\pb{s^*H,t^*C}=0$, so these vector fields commute; moreover, they
are independent at a a generic point of ${\rm SO}(3) \times {\bbR}^3$. On the other hand, the functions $t^*x,t^*y$
and $t^*z$ are in involution with $s^*H$ as well as with $t^*C$. It follows that\footnote{In this list of functions
  one can replace $t^*y$ or $t^*z$ by $t^*x$.} $(s^*H,t^*C,t^*y,t^*z)$ defines a non-commutative integrable system
of rank 2 on ${\rm SO}(3) \times {\bbR}^3$.

For our purposes we need to restrict phase space to an open subset on which the NCI system is regular. Let us
denote by $\vert\vert\cdot\vert\vert$ the standard norm on $\bbR^3$, so for $m=(x,y,z)\in\bbR^3$ we have
$\vert\vert m\vert\vert^2=x^2+y^2+z^2$. The inequalities $I_x>I_y>I_z>0$ imply that the image of $H$ is the closed
interval
\begin{equation*}
  \Im(\phi_H)=\set{(v,h)\mid\frac{\vert\vert m\vert\vert^2}{2I_x}\leqs h\leqs\frac{\vert\vert m\vert\vert^2}{2I_z}}\;.
\end{equation*}%
Let $B$ and $B'$ denote the open subsets of $\bbR^3\times \bbR$, defined by
\begin{eqnarray*}
  B&:=&\left\{(v,h)\mid\frac{\vert\vert v\vert\vert^2}{2I_x}<h<\frac{\vert\vert v\vert\vert^2}{2I_y}\right\}\;, \\
  B'&:=&\left\{(v,h)\mid\frac{\vert\vert v\vert\vert^2}{2I_y}<h<\frac{\vert\vert v\vert\vert^2}{2I_z}\right\}\;.
\end{eqnarray*}%
We denote by $M\subset {\rm SO}(3) \times {\bbR}^3$ the inverse image $\phi_H^{-1}(B)$, consisting of all $(R,m)$
for which $(m,H(m))\in B$; the analysis done below can be repeated with minor changes for $M':=\phi_H^{-1}(B')$.
On $M$ the NCI system is regular; more precisely $(M,\Pi)\stackrel{\phi_H}{\to}(B,\pi)$ is a rank two NCI system
with momentum map. The fibers of $\phi_H$ are compact but not connected: the fiber over each point of $B$ consists
of two disjoint two-dimensional tori $\bbT^2$. Since, for our analysis, we need the fibers of the momentum map to
be connected, we need to do a further restriction on phase space: we define $M_+$ as the subset of $M$ whose points
$(R,m)$, with $m=(x,y,z)$, satisfy $x> 0$. Now $(M_+,\Pi)\stackrel{\phi_H}{\to}(B,\pi)$ is a regular NCI system
of rank two with compact connected fibers.

For explicitness, we give a geometrical description of these fibers as two-dimensional tori. Let $(v,h)\in
B\subset\bbR^3\times\bbR$ and let $c:=\vert\vert v\vert\vert$. The fiber in $M_+$ over $(v,h)$ is given by
\begin{equation*}
  \phi_H^{-1}(v,h)=\set{(R,m)\in {\rm SO}(3) \times {\bbR}^3\mid Rm=v,\,H(m)=h}\;.
\end{equation*}%
Notice that when $(R,m)\in \phi_H^{-1}(v,h)$, the point $m$ belongs to one of the two connected components of the
intersection of the sphere $\vert\vert m\vert\vert^2=c$ and the ellipsoid $H(m)=h$. This component, which
corresponds to the component lying in the half-space $x>0$ (see the above definition of $M_+$) is a smooth curve
$S$, diffeomorphic\footnote{The \emph{complex} intersection of these two quadrics is a smooth complex elliptic
  curve.}  to the circle $S^1$. Notice also that if $R_v$ is any rotation with center $O$ which fixes $v$ then
$(R_vR,m)$ belongs to the same fiber of $\phi_H$. This leads to two actions of $S^1$ on $\phi_H^{-1}(v,h)$. The
first one leaves $m$ unchanged and is the above left multiplication of $R$ by the unique rotation $R_v$ over a
given angle. For the action of the other component $S^1$ one fixes a diffeomorphism between $S$ and $S^1$; the
action on $m$, denoted $\theta\cdot m$ is then given by the standard action of $S^1$ on itself, while the action on
$R$ can be taken as right multiplication of $R$ with the unique rotation which sends $\theta\cdot m$ to
$m$. Clearly these two actions of $S^1$ commute and they define an action of $\bbT^2$ which is transitive and has
trivial stabilizer. It allows us to identify (topologically) $\phi_H^{-1}(v,h)$ with $\bbT^2$.

We now address the question of the existence of action-angle variables and foliations for the Euler-Poinsot top (on
$M_+$). First, since $M_+$ is a symplectic manifold, the symplectic foliation on $B$ is regular and is the only
action foliation (see Remark \ref{rem:action_foliation_symplectic}), in particular there exists an action
foliation. Moreover, since $B$ is simply-connected there are no obstructions to extend the action variables which
define locally the action foliation into global action variables. Thus, global action variables exist also.

We finally show that the Euler-Poinsot system does not admit an angle foliation, hence does not admit global angle
variables. To do this, we show that the submersion $\phi_H:M_+\to B$ does not admit a coisotropic section. Notice
first that $B$ is, topologically, the product of a $2$-sphere by $\bbR$. In particular, it is simply-connected,
i.e. $\pi_1(B)=0$, but it is not $2$-connected, i.e. $\pi_2(B)$ is not trivial. On the contrary, 
$$ 
  M_+=SO(3)\times\set{(x,y,z)\neq(x,0,0)\mid x>0\hbox{ and } x^2<\frac{I_y-I_z}{I_x-I_y}\frac{I_x}{I_z}z^2}
$$
from which we see that $M_+$ is homeomorphic to $SO(3)\times\bbR_{>0}\times(\bbR^2\setminus\set0)$, so that $M_+$
is $2$-connected but not simply-connected.  The argument is now purely topological. Assume that an angle foliation
exists, and denote by $ {\mathcal F}$ one of its leaves.  By construction, ${\mathcal F} $ is a connected
submanifold and the restriction of $ \phi_H$ to $ {\mathcal F}$ is a local diffeomorphism onto $B$. Since $B$ is
simply-connected, the restriction of $ \phi_H$ to $ {\mathcal F}$ has to be a global diffeomorphism. Inverting the
restriction of $ \phi_H$ to ${\mathcal F} $ yields a global section of $ \phi_H$. But this is in turn impossible
because $\pi_2(B)$ is not trivial while $ \pi_2(M_+)$ is trivial, which prohibits the existence of such a
section. Hence the Euler-Poinsot top admits neither a set of angle variables nor an angle foliation.The fact that
angle variables for the Euler-Poinsot do not exist was already shown by F. Fasso (see \cite{FA}).

\subsection{The Gelfand-Cetlin system}
We finish with a non-trivial example where action-angle variables exist: the Gelfand-Cetlin system. The results in
this section are due to A. Giacobbe and we refer to his orignal paper \cite{Giacobbe} for details and proofs.

The phase space of the Gelfand-Cetlin system is the real vector space of $n\times n$ hermitian matrices $\fH_n$. It
has a linear Poisson structure, since it can be viewed as the dual of the Lie algebra of unitary matrices
$\fu_n$. Explicitly, the Poisson structure $\Pi$ is given for smooth functions $F,G$ on $\fH_n$ at $X\in\fH_n$ by
\begin{equation*}%
  \pb{F,G}(X):=\inn{\lb{\nabla F(X),\nabla G(X)}}{X}\;,  
\end{equation*}%
where the inner product is defined for $X,Y\in\fH_n$ by $\inn XY:=i\Tr XY$ and $\nabla F(X)$ is the differential of
$F$ at $X$, viewed as an element of $\fH_n$ (using the inner product). The rank of this Poisson structure is
$n(n-1)$, to be compared with $\dim\fH_n=n^2$. When one removes from $X\in\fH_n$ the last $n-i$ rows and columns
one obtains an element of $\fH_{n-i}$, which is denoted by $X^{(i)}$. For $i=1,\dots,n$ the $i$ eigenvalues of
$X^{(i)}$ are denoted by $\mu_p^i(X)$; they are ordered such that $\mu_1^i(X)\leqs
\mu_2^i(X)<\cdots\leqs\mu_{n-i}^i(X)$. They satisfy
\begin{equation}\label{eq:GC_ineqs}
  \mu_p^{i+1}(X)\leqs\mu_p^i(X)\leqs\mu_{p+1}^{i+1}(X)\;.
\end{equation}%
Let $M$ be the open subset of $\fH_n$ where each $X^{(i)}$ has simple spectrum and where the eigenvalues of
$X^{(i)}$ are different from the eigenvalues of $X^{(i+1)}$. On~$M$ the maps $X\mapsto \mu_p^i(X)$ define
$N:=n(n+1)/2$ smooth functions, which are independent, leading to a submersion $\phi:M\to B$, where $B$ is the
sector in $\bbR^N$, defined by replacing in (\ref{eq:GC_ineqs}) the inequalities by strict inequalities. Moreover,
these functions are in involution and the NCI system $ (M,\Pi) \stackrel{\phi}{\to} (B ,0)$ is regular. The fibers
of $\phi$ are compact and connected, i.e., they are diffeomorphic to tori of dimension $r:=n(n-1)/2=\Rk\Pi/2$.

The $n$ functions $\mu_1^n,\mu_2^n,\dots,\mu_n^n$ are Casimirs of $\Pi$, while the other $N/2$ functions $\mu_p^i$
($i<n$) have independent periodic flows of period 1. Thus, they provide a set of action variables. The construction
of the angles variables is slightly more involved. For given $i$ such that $0<i<n$ we explain how to compute the
angle variables $\varphi_p^i$ which are conjugate to $\mu_p^i$, for $p=1,\dots,p$. The main operation involved in
computing $\varphi_p^i(X)$ for $X\in\fH_n$ is to conjugate~$X$ by a unitary block matrix of the form
$\Lambda:=\left(\begin{matrix}P&0\\0&I_{n-i}\end{matrix}\right)$ such that $\Lambda X\bar\Lambda^t$ is of the form
$X':=\left(\begin{matrix}\Delta&*\\ *&*\end{matrix}\right)$, where $\Delta$ is diagonal, i.e., $\Delta
=\mathop{\rm diag}(\mu_1^i,\dots,\mu_i^i)$. Of course, such a matrix $P$ is not unique, but all entries of its last row are
non-zero and a unique $P$ can be selected by demanding that all these entries are strictly positive real numbers
and that the columns have norm 1. With this choice of $P$, the angle variable $\varphi_p^i(X)$ is the argument of
the complex number $X'_{p,i+1}$. Combined, the set of $(\mu_p^i,\varphi_p^i)$, where $i$ ranges from~$1$ to $n-1$
and $p$ from~$1$ to $i$, provide a set of action-angle variables for the Gelfand-Cetlin system.
\bibliographystyle{abbrv}
\bibliography{ref}

\def\cprime{$'$}
\begin{thebibliography}{10}

\bibitem{adlermoerbekevanhaecke2004}
M.~Adler, P.~van Moerbeke, and P.~Vanhaecke.
\newblock {\em Algebraic integrability, {P}ainlev\'e geometry and {L}ie
  algebras}, volume~47 of {\em Ergebnisse der Mathematik und ihrer
  Grenzgebiete. 3. Folge. A Series of Modern Surveys in Mathematics [Results in
  Mathematics and Related Areas. 3rd Series. A Series of Modern Surveys in
  Mathematics]}.
\newblock Springer-Verlag, Berlin, 2004.

\bibitem{arnold}
V.~Arnol{\cprime}d.
\newblock {\em Mathematical methods of classical mechanics}.
\newblock Springer-Verlag, New York, 1978.
\newblock Translated from the Russian by K. Vogtmann and A. Weinstein, Graduate
  Texts in Mathematics, 60.

\bibitem{bolsinov}
A.~V. Bolsinov and B.~Jovanovi{\'c}.
\newblock Noncommutative integrability, moment map and geodesic flows.
\newblock {\em Ann. Global Anal. Geom.}, 23(4):305--322, 2003.

\bibitem{candel}
A.~Candel and L.~Conlon.
\newblock {\em Foliations. {I}}, volume~23 of {\em Graduate Studies in
  Mathematics}.
\newblock American Mathematical Society, Providence, RI, 2000.

\bibitem{Cannas_weinstein}
A.~Cannas~da Silva and A.~Weinstein.
\newblock {\em Geometric models for noncommutative algebras}, volume~10 of {\em
  Berkeley Mathematics Lecture Notes}.
\newblock American Mathematical Society, Providence, RI; Berkeley Center for
  Pure and Applied Mathematics, Berkeley, CA, 1999.

\bibitem{CDW}
A.~Coste, P.~Dazord, and A.~Weinstein.
\newblock Groupo\"\i des symplectiques.
\newblock In {\em Publications du {D}\'epartement de {M}ath\'ematiques.
  {N}ouvelle {S}\'erie. {A}, {V}ol.\ 2}, volume~87 of {\em Publ. D\'ep. Math.
  Nouvelle S\'er. A}, pages i--ii, 1--62. Univ. Claude-Bernard, Lyon, 1987.

\bibitem{batesandcushman}
R.~H. Cushman and L.~M. Bates.
\newblock {\em Global aspects of classical integrable systems}.
\newblock Birkh\"auser Verlag, Basel, 1997.

\bibitem{dazorddelzant}
P.~Dazord and T.~Delzant.
\newblock Le probl\`eme g\'en\'eral des variables actions-angles.
\newblock {\em J. Differential Geom.}, 26(2):223--251, 1987.

\bibitem{MR1104924}
J.-P. Dufour and P.~Molino.
\newblock Compactification d'actions de {${\bf R}^n$} et variables action-angle
  avec singularit\'es.
\newblock In {\em Symplectic geometry, groupoids, and integrable systems
  ({B}erkeley, {CA}, 1989)}, volume~20 of {\em Math. Sci. Res. Inst. Publ.},
  pages 151--167. Springer, New York, 1991.

\bibitem{Duis}
J.~Duistermaat.
\newblock On global action-angle coordinates.
\newblock {\em Comm. Pure Appl. Math.}, 33(6):687--706, 1980.

\bibitem{FA}
F.~Fass{\`o}.
\newblock The {E}uler-{P}oinsot top: a non-commutatively integrable system
  without global action-angle coordinates.
\newblock {\em Z. Angew. Math. Phys.}, 47(6):953--976, 1996.

\bibitem{AF}
E.~Fiorani and G.~Sardanashvily.
\newblock Noncommutative integrability on noncompact invariant manifolds.
\newblock In {\em X{V} {I}nternational {W}orkshop on {G}eometry and {P}hysics},
  volume~11 of {\em Publ. R. Soc. Mat. Esp.}, pages 282--286. R. Soc. Mat.
  Esp., Madrid, 2007.

\bibitem{Giacobbe}
A.~Giacobbe.
\newblock Some remarks on the {G}elfand-{C}etlin system.
\newblock {\em J. Phys. A}, 35(49):10591--10605, 2002.

\bibitem{MR1725760}
J.~Grabowski, G.~Marmo, and P.~W. Michor.
\newblock Construction of completely integrable systems by {P}oisson mappings.
\newblock {\em Modern Phys. Lett. A}, 14(30):2109--2118, 1999.

\bibitem{guilleminandsternberg}
V.~Guillemin and S.~Sternberg.
\newblock The {G}el\cprime fand-{C}etlin system and quantization of the complex
  flag manifolds.
\newblock {\em J. Funct. Anal.}, 52(1):106--128, 1983.

\bibitem{LMV}
C.~Laurent-Gengoux, E.~Miranda, and P.~Vanhaecke.
\newblock Action-angle coordinates for integrable systems on {P}oisson
  manifolds.
\newblock {\em Int. Math. Res. Not. IMRN}, (8):1839--1869, 2011.

\bibitem{liber}
P.~Libermann.
\newblock Probl\`emes d'\'equivalence et g\'eom\'etrie symplectique.
\newblock In {\em Third {S}chnepfenried geometry conference, {V}ol. 1
  ({S}chnepfenried, 1982)}, volume 107 of {\em Ast\'erisque}, pages 43--68.
  Soc. Math. France, Paris, 1983.

\bibitem{MR882548}
P.~Libermann and C.-M. Marle.
\newblock {\em Symplectic geometry and analytical mechanics}, volume~35 of {\em
  Mathematics and its Applications}.
\newblock D. Reidel Publishing Co., Dordrecht, 1987.
\newblock Translated from the French by Bertram Eugene Schwarzbach.

\bibitem{Nekhroshev}
N.~N. Nehoro{\v{s}}ev.
\newblock Action-angle variables, and their generalizations.
\newblock {\em Trudy Moskov. Mat. Ob\v s\v c.}, 26:181--198, 1972.

\bibitem{sussmann}
H.~J. Sussmann.
\newblock Orbits of families of vector fields and integrability of
  distributions.
\newblock {\em Trans. Amer. Math. Soc.}, 180:171--188, 1973.

\bibitem{whittaker}
E.~T. Whittaker.
\newblock {\em A treatise on the analytical dynamics of particles and rigid
  bodies}.
\newblock Cambridge Mathematical Library. Cambridge University Press,
  Cambridge, 1988.
\newblock With an introduction to the problem of three bodies, Reprint of the
  1937 edition, With a foreword by William McCrea.

\bibitem{MR1389366}
N.~T. Zung.
\newblock Symplectic topology of integrable {H}amiltonian systems. {I}.
  {A}rnold-{L}iouville with singularities.
\newblock {\em Compositio Math.}, 101(2):179--215, 1996.

\end{thebibliography}

\end{document}